\documentclass[twocolumn]{autart}

\usepackage{dsfont}
\usepackage{amsfonts}
\usepackage{mathrsfs}
\usepackage{}
\usepackage{colortbl}
\usepackage{amssymb,amsmath,graphicx}
\usepackage{float}
\usepackage{multirow}

\newtheorem{Remark}{Remark}[section]

\newtheorem{Problem}{Problem}
\newenvironment{Proof}{\noindent{\em Proof:\/}}{\hfill $\Box$\par}
\newtheorem{Theorem}{Theorem}[section]
\newtheorem{Proposition}{Proposition}[section]
\newtheorem{Lemma}{Lemma}[section]

\newtheorem{Assumption}{Assumption}

\newcommand{\EQQ}{\begin{eqnarray*}}
\newcommand{\ENN}{\end{eqnarray*}}
\newcommand{\EQ}{\begin{eqnarray}}
\newcommand{\EN}{\end{eqnarray}}

\begin{document}

\begin{frontmatter}

\title{Cooperative Global Robust Output Regulation for a Class of Nonlinear Multi-Agent Systems by Distributed Event-Triggered Control
\thanksref{footnoteinfo}} 

\thanks[footnoteinfo]{This work has been supported  by the Research Grants Council of the Hong Kong Special Administration Region under grant
No. 14200515.
Corresponding author: Jie Huang.
}

\author[Liu-Huang]{Wei Liu}\ead{wliu@mae.cuhk.edu.hk},
\author[Liu-Huang]{Jie Huang}\ead{jhuang@mae.cuhk.edu.hk}

\address[Liu-Huang]{Department of Mechanical and Automation Engineering, The Chinese University of Hong Kong, Shatin, N.T., Hong Kong}


\begin{abstract}
This paper studies the event-triggered cooperative global robust output regulation problem for a class of nonlinear multi-agent systems via a distributed internal model design.
We show that our problem can be solved practically in the sense
that the ultimate bound of the tracking error can be made arbitrarily small by adjusting a design parameter in the proposed event-triggered mechanism. Our result offers a few new features. First, our control law is robust against both external disturbances and parameter uncertainties, which are allowed to belong to some arbitrarily large prescribed compact sets.
Second,  the nonlinear functions in our system do not need to satisfy the global Lipchitz condition. Thus our systems are general enough to include some benchmark nonlinear systems that cannot be handled by existing approaches. Finally, our control law is a specific distributed output-based event-triggered control law, which lends itself to a direct digital implementation.
\end{abstract}

\end{frontmatter}

\section{Introduction}
Over the past decade, various cooperative control problems for multi-agent systems have been widely studied. The cooperative robust output regulation problem is one of the fundamental and important cooperative control problems, which aims to make all followers track some reference inputs and reject some external disturbances, where the  reference inputs and the external disturbances are both generated by an exosystem called the leader system. The problem has been first studied for linear uncertain multi-agent systems in \cite{SuHongHuang2013,Wang3}, and then for nonlinear uncertain multi-agent systems in \cite{Ding1,DongHuang2014a,SuHuang2015}. In this paper,  we will further study  the event-triggered cooperative global robust practical output regulation problem for a class of nonlinear multi-agent systems in normal form with unity relative degree. 

Our study is motivated by the need of implementing  continuous-time control laws in digital platforms. Compared with the traditional sampled-data implementation \cite{Astrom1,Franklin1}, where the data-sampling is performed periodically, the  event-triggered control approach generates the samplings and control actuation depending on the real time system state or output, and is more efficient in reducing the number of control task executions while maintaining the control performance \cite{Heemels1}.
So far, extensive efforts have been made for event-triggered control for single linear systems in, e.g.,  \cite{Donkers2012,Heemels1} and for nonlinear systems in, say, \cite{Girard1,LiuT1,LiuHuang2017a,Tabuada1,XingL1}. In particular,
reference \cite{Heemels1} gave an introduction to the event-triggered control and studied the stabilization problem for a class of linear systems by a state-feedback event-triggered control law.
Reference \cite{Donkers2012} analyzed the  closed-loop  stability and the $\mathcal{L}_{\infty}$-performance for a class of linear systems by an output-based event-triggered control law.
In \cite{Tabuada1},  the  stabilization problem for a class of nonlinear systems was solved by a state-based event-triggered control law.
Reference \cite{Girard1} further proposed a dynamic event-triggered mechanism to solve the stabilization problem for the same class of nonlinear systems as that in \cite{Tabuada1}.
In \cite{LiuT1},  a state-based event-triggered control law was designed to solve the robust stabilization problem for a class of nonlinear systems subject to external disturbances by applying  the cyclic small gain theorem.
Reference \cite{XingL1} solved the tracking problem for a class of high-order uncertain nonlinear systems  by  an event-triggered adaptive control law.
In \cite{LiuHuang2017a},   the global robust output regulation problem for a class of nonlinear systems  was further solved by an  output-based event triggered control law.
Other relevant contributions can also be found in \cite{Abdelrahim1,Dolk2016,Postoyan1,Tallapragada1,XingL1} etc.

The event-triggered control approach has also been applied to the cooperative control problems for multi-agent systems. For example, the consensus problem was studied  by the event-triggered control approach for linear multi-agent systems in \cite{Cheng1,Dimarogonas1,Fan1,Seyboth1,Zhu1}. 
In \cite{HuW2},  based on the feedforward design, the cooperative output regulation problem for a class of exactly known linear multi-agent systems was studied by a distributed  event-triggered control law. Reference \cite{LiuHuang2017b} further designed a distributed output-based event-triggered control law to solve the cooperative robust output regulation problem for a class of minimum phase linear uncertain multi-agent systems  based on the internal model approach. Also, references \cite{LiH2,LiH3,XieD1} studied the event-triggered consensus problem  for several types of nonlinear multi-agent systems satisfying the global Lipchitz condition.

Compared with the existing results on the event-triggered cooperative control problems for nonlinear multi-agent systems in \cite{LiH2,LiH3,XieD1}, our problem offers at least three new features. First, our system contains both external disturbances and parameter uncertainties, which are allowed to belong to some arbitrarily large prescribed compact sets.
Second,  the nonlinear functions in our system do not need to satisfy the global Lipchitz condition. Thus our systems are general enough to include some benchmark nonlinear systems such as Lorenz systems, FitzHugh-Nagumo systems, which cannot be handled by
the approaches in \cite{LiH2,LiH3,XieD1}. Finally, our control law is a dynamic distributed output-based event-triggered control law, which is more challenging than the static distributed state-based or static distributed output-based control law, since we need to sample not only the measurement output of each agent but also the state of the dynamic compensator.   To overcome these challenges, we combine the distributed internal model approach with a distributed event-triggered mechanism. This event-triggered mechanism contains a design parameter
which not only dictates the ultimate bound of the closed-loop system, but also the frequency of the triggering events.
It is shown that this control law together with the event-triggered mechanism solves our problem
in the sense that the steady-state tracking error of the closed-loop system can be made arbitrarily small. Besides, our method guarantees the existence of the minimal inter-execution time of the event-triggered mechanism, thus preventing the Zeno behavior from happening.

It is worth mentioning that the problem in this paper also contain the problem in \cite{LiuHuang2017a} as a special case by letting $N=1$. Compared with \cite{LiuHuang2017a}, the main challenge of this paper is that we need to design a specific event-triggered control law and a specific event-triggered mechanism satisfying some communication constraints for each subsystem, where the communication constraints means that each subsystem can only make use of the information of its neighbors and itself for control. We call an  event-triggered control law and an event-triggered mechanism satisfying such communication constraints as a distributed event-triggered control law and a distributed event-triggered mechanism, respectively.

{\bf Notation.} For any column vectors $a_i$, $i=1,...,s$, denote $\mbox{col}(a_1,...,a_s)=[a_1^T,...,a_s^T]^T$.
$\|x\|$ denotes the Euclidean norm of vector $x$. $\|A\|$ denotes the induced norm of matrix $A$ by the Euclidean norm. $\mathbb{Z}^+$ denotes
the set of nonnegative integers.
  $\lambda_{\max}(A)$ and $\lambda_{\min}(A)$ denote the maximum eigenvalue and the minimum eigenvalue of a symmetric real matrix $A$, respectively.
A matrix $M\in\mathbb{R}^{N\times N}$ is called an $\mathcal{M}-$matrix if all of its non-diagonal elements are non-positive and all of its eigenvalues have positive real parts.

\section{Problem formulation and preliminaries}\label{PF}
Consider a class of nonlinear multi-agent systems taken from \cite{DongHuang2014a} as follows:
\begin{equation}\label{system1}
\begin{split}
  \dot{z}_{i} =& f_{i}(z_{i},y_{i},v,w)\\
  \dot{y}_{i} =& g_{i}(z_{i},y_{i},v,w) + b_{i}(w)u_{i}\\
  e_{i}=&y_{i}-q(v,w),\ \ \ i=1,\cdots,N
\end{split}
\end{equation}
where, for $i=1,\cdots,N$, $(z_{i},y_{i})\in \mathbb{R}^{n}\times\mathbb{R}$ is the state,  $e_{i}\in \mathbb{R}$ is the error output, $u_{i}\in \mathbb{R}$ is the input,
$w\in\mathbb{R}^{n_{w}}$ is an uncertain constant vector, and $v(t)\in\mathbb{R}^{n_{v}}$ is an exogenous signal representing both reference input to be tracked and disturbance to be rejected.
Here $v(t)$ is assumed to be generated by  the following linear system:
\begin{equation}\label{exosystem1}
\begin{split}
  \dot{v} = Sv,\ \ y_{0}=q(v,w).
\end{split}
\end{equation}
We assume that all functions in (\ref{system1}) and (\ref{exosystem1}) are sufficiently smooth, and satisfy $b_{i}(w)>0$, $f_{i}(0,0,0,w)=0$, $g_{i}(0,0,0,w)=0$ and $q(0,w)=0$ for all
$w\in \mathbb{R}^{n_{w}}$.

System \eqref{system1} is the so-called  nonlinear multi-agent system in normal form with unity relative degree.  Like in \cite{DongHuang2014a}, the plant (\ref{system1}) and the
exosystem (\ref{exosystem1}) together can be viewed as a multi-agent system of $N+1$ agents with (\ref{exosystem1}) as the leader and the $N$
subsystems of (\ref{system1}) as $N$ followers.
Given the plant (\ref{system1}) and the exosystem (\ref{exosystem1}),
  we can define a digraph $\bar{\mathcal{G}}=(\bar{\mathcal{V}},\bar{\mathcal{E}})$
where $\bar{\mathcal{V}}=\{0,1,\cdots,N\}$ with $0$ associated with
the leader system and with $ i = 1,\cdots,N$ associated with the $N$
followers, respectively, and $\bar{\mathcal{E}}\subseteq
\bar{\mathcal{V}}\times\bar{\mathcal{V}}$. For each $j=0,1,\cdots,N$, $i=1,\cdots,N$, and $i \neq j$,
$(j,i)\in\bar{\mathcal{E}} $  if and only if the control
$u_{i}$ can make use of $y_{i}-y_{j}$ for feedback
control. If the digraph $\bar{\mathcal {G}} $ contains a sequence of edges of the form $(i_{1},i_{2}),(i_{2},i_{3}),\cdots,(i_{k},i_{k+1})$, then the node $i_{k+1}$ is said to be reachable from the node $i_{1}$. For $i=1,\cdots,N$, let $\bar{\mathcal{N}}_{i}=\{j, (j,i)\in\bar{\mathcal{E}}\}$ denote the neighbor set of node $i$.

Denote the adjacency matrix of the digraph
$\bar{\mathcal{G}}$ by $\bar{\mathcal{A}} =[\bar{a}_{ij} ]\in \mathbb{R}^{(N+1)\times
(N+1)}$ where $\bar{a}_{ii} =0$, $\bar{a}_{ij} =1\Leftrightarrow
(j,i)\in\mathcal{\bar{E}}$, and $\bar{a}_{ij} =0\Leftrightarrow
(j,i)\notin\mathcal{\bar{E}}$ for $i,j=0,1,\cdots,N$.
Define the virtual error output for agent $i$ as
\begin{equation}\label{evi1}
\begin{split}
e_{vi}(t)=\sum_{j=0}^{N}\bar{a}_{ij}(y_{i}(t)-y_{j}(t)),\ i=1,\cdots,N.
\end{split}
\end{equation}
Let $e_{0}=0$, $e=\mbox{col} (e_{1},\cdots,e_{N})$, $e_{v}=\mbox{col} (e_{v1},\cdots,e_{vN})$, and
$H=[h_{ij}]_{i,j=1}^{N}$ with $h_{ii} =\sum_{j=0}^{N}\bar{a}_{ij}$ and $h_{ij}=-\bar{a}_{ij}$ for $i\neq j$.
It can be easily verified that $e_{v}=He$.

For $i=1,\cdots,N$,
consider a control law of the following form
\begin{equation}\label{ui1}
\begin{split}
  &u_{i}(t)= \hat{f}_{i}(\eta_{i}(t_{k}^{i}),e_{vi}(t_{k}^{i}))\\
  &\dot{\eta}_{i}(t)= \hat{g}_{i}(\eta_{i}(t),e_{vi}(t_{k}^{i})),~\forall t\in[t_{k}^{i},t_{k+1}^{i}),~k\in\mathbb{Z}^{+}
\end{split}
\end{equation}
where   $\hat{f}_{i}(\cdot)$ and $\hat{g}_i(\cdot)$ are some nonlinear functions, $\eta_{i}$ subsystem is the so-called internal model which will be
designed later, $t_{k}^{i}$ denotes the triggering time instants of agent $i$ with $t_{0}^{i}=0$, and is generated by the following event-triggered mechanism
\begin{equation}\label{trigger1}
\begin{split}
t_{k+1}^{i}=\inf\{t>t_{k}^{i}~|~\hat{h}_{i}(\tilde{e}_{vi}(t),\tilde{\eta}_{i}(t),e_{vi}(t))\geq\delta\}
\end{split}
\end{equation}
where  $\hat{h}_{i}(\cdot)$ is some nonlinear function, $\delta>0$ is some constant, and
 \begin{equation}\label{tildeevi1}
\begin{split}
 &\tilde{e}_{vi}(t)=e_{vi}(t_{k}^{i})-e_{vi}(t)\\
 &\tilde{\eta}_{i}(t)=\eta_{i}(t_{k}^{i})-\eta_{i}(t)\\
\end{split}
\end{equation}
for any $t\in[t_{k}^{i},t_{k+1}^{i})$ with $i=1,\cdots,N$ and $k\in\mathbb{Z}^{+}$.


\begin{Remark}
 It is noted both the control law \eqref{ui1} and the event-triggered mechanism \eqref{trigger1} satisfy the communication constraints in the sense
 that the control law and the event-triggered mechanism of each agent only make use of the output information of its neighbors and itself. Thus,
We call a control law of the form \eqref{ui1} and an event-triggered mechanism of the form \eqref{trigger1} as a distributed output-based event-triggered control law and  a distributed output-based event-triggered mechanism, respectively.
It can be seen that our control law and event-triggered mechanism have the following three advantages over those in the existing references.
First, the distributed event-triggered mechanism \eqref{trigger1} is more practical than the centralized event-triggered mechanism used in \cite{XieD1}, where the event-triggered mechanism of each agent depends on the full state information of all agents. Second, since the event-triggered mechanism \eqref{trigger1} works in an asynchronous way, that is, the triggering time instants of each agent are  generated independently of the triggering time instants of all other agents, it is more efficient than the event-triggered mechanism used in \cite{Dimarogonas1} which works in a synchronous way, i.e., the control laws of all agents are updated simultaneously. Third, since our control law takes a specific form and generates piecewise constant signals, as will be seen in Remark \ref{Discretization}, it  can be directly implemented in a digital platform.
\end{Remark}
\begin{Remark}
The event-triggered mechanism (\ref{trigger1}) is said to have a minimal inter-execution time if there exists a real number $\tau_{d} >0$ such that
$t_{k+1}^{i} - t_{k}^{i} \geq \tau_{d}$ for all $k \geq 0$ and $i=1,\cdots,N$.
 Clearly,  if such a minimal inter-execution time $\tau_{d}$ exists, then  the execution times of the control law cannot become arbitrarily close. Thus,  the Zeno behavior can be avoided.
\end{Remark}

Now we describe our problems as follows:
%
\begin{Problem}\label{Problem1}
 Given the plant (\ref{system1}), the exosystem (\ref{exosystem1}),  a  digraph $\bar{\mathcal{G}}$,  some compact subsets
$\mathbb{V}\subset\mathbb{R}^{n_{v}}$ and $\mathbb{W}\subset\mathbb{R}^{n_{w}}$ with $0\in\mathbb{W}$ and $0\in\mathbb{V}$, and any $\epsilon>0$, design an event-triggered mechanism of the form \eqref{trigger1}  and a control law of the form (\ref{ui1})  such that the resulting closed-loop system has the following properties: for any $v\in\mathbb{V}$, $w\in\mathbb{W}$, and any initial states
$z_{i}(0)$, $y_{i}(0)$, $\eta_{i}(0)$,
\begin{enumerate}
  \item  the trajectory of the closed-loop system  exists and is bounded for all $t\geq0$;
  \item $\lim_{t \to \infty}\sup\|e(t)\|\leq \epsilon$.
\end{enumerate}
\end{Problem}

\begin{Remark}\label{RemarkProb1}
We call Problem \ref{Problem1}  the event-triggered cooperative global robust practical output regulation problem
and call a control law that solves Problem \ref{Problem1}  a practical solution to the cooperative global robust output regulation problem.
Such a problem offers
at least four new features compared with the problem in studied in \cite{DongHuang2014a}, where the cooperative output regulation problem
for the  nonlinear multi-agent system was studied by an analog control law. First, here we  need to design not only a piecewise constant control law but
also an event-triggered mechanism, while, in \cite{DongHuang2014a}, only a continuous-time control law needs to be designed. Second, since, under the event-triggered control, the closed-loop
system is a hybrid system, the stability analysis of the closed-loop system is much more sophisticated than the one in \cite{DongHuang2014a}.   Third, the Zeno
phenomenon is unique for the closed-loop system under the event-triggered control and our control law and event-triggered mechanism need to be designed to exclude the Zeno phenomenon. Finally, as will be seen in next section,  the topological assumption will be relaxed from the undirected connected network case to the more general directed connected network case, which further complicated the stability analysis of the closed-loop system.
\end{Remark}

It is known from the framework in \cite{Huang2} that the output regulation problem for a given plant can be converted to a
stabilization problem of a well defined augmented system.  In order to form the so-called
augmented system,
we first introduce some standard assumptions which can also be found in \cite{DongHuang2014a}.

\begin{Assumption}\label{Ass2.1}
The exosystem is neutrally stable, i.e., all the eigenvalues of S are semi-simple with zero real parts.
\end{Assumption}

\begin{Remark}
Assumption \ref{Ass2.1} has also been  used in \cite{Ding1,DongHuang2014a,SuHuang2015}. Under Assumption \ref{Ass2.1}, for any $v(0)\in\mathbb{V}_{0}$ with $\mathbb{V}_{0}$ being a compact set, there exists another compact set $\mathbb{V}$ such that $v(t)\in\mathbb{V}$ for all $t\geq0$.
\end{Remark}

%

\begin{Assumption}\label{Ass2.3}
There exist globally defined smooth functions $\textbf{z}_{i}:\mathbb{R}^{n_{v}}\times\mathbb{R}^{n_{w}}\mapsto\mathbb{R}^{n}$ with $\textbf{z}_{i}(0,w)=0$ such that
\begin{equation}\label{Zi}
\begin{split}
  \dfrac{\partial\textbf{z}_{i}(v,w)}{\partial v}Sv=f_{i}(\textbf{z}_{i}(v,w),q(v,w),v,w)
\end{split}
\end{equation}
for all $(v,w)\in\mathbb{R}^{n_{v}}\times\mathbb{R}^{n_{w}}$, $i=1,\cdots,N$.
\end{Assumption}
\begin{Remark}\label{RemarkAss2.3}
Under Assumption \ref{Ass2.3}, let
\begin{equation}\label{Ui}
\begin{split}
\textbf{y}_{i}(v,w)=&~q(v,w)\\
  \textbf{u}_{i}(v,w)=&~b_{i}^{-1}(w)\big(\dfrac{\partial q(v,w)}{\partial v}Sv\\
  &-g_{i}(\textbf{z}_{i}(v,w),q(v,w),v,w)\big).
\end{split}
\end{equation}
Then, $\textbf{z}_{i}(v,w)$, $\textbf{y}_{i}(v,w)$ and $\textbf{u}_{i}(v,w)$ are the solutions to the regulator equations
associated with (\ref{system1}) and (\ref{exosystem1}).  It can be seen that Assumption \ref{Ass2.3} is a necessary condition for the solvability of the regulator equations
associated with (\ref{system1}) and (\ref{exosystem1}) and thus a necessary condition for the solvability of the cooperative output regulation problem of \eqref{system1} and \eqref{exosystem1}  \cite{Huang1}.
\end{Remark}

\begin{Assumption}\label{Ass2.4}
 $\textbf{u}_{i}(v,w),\ i=1,\cdots,N$, are polynomials in $v$ with coefficients depending on $w$.
\end{Assumption}

 \begin{Remark}\label{RemarkAss2.4}
As remarked in \cite{DongHuang2014a}, under Assumptions \ref{Ass2.1} to \ref{Ass2.4}, for $i=1,\cdots,N$, there exist some integers $s_{i}$
and some real coefficients polynomials
\begin{equation}\label{bPi1}
\begin{split}
P_{i}(\lambda)=\lambda^{s_{i}}-\varrho_{1i}-\varrho_{2i}\lambda-\cdots-\varrho_{s_{i}i}\lambda^{s_{i}-1}
\end{split}
\end{equation}
 whose roots are all distinct with zero real part,
 such that, for all trajectories $v(t)$ of the exosystem and all $w\in\mathbb{W}$, $\textbf{u}_{i}(v,w)$ satisfy
\begin{equation}\label{dotUi}
\begin{split}
   \dfrac{d^{s_{i}}\textbf{u}_{i}}{dt^{s_{i}}}=\varrho_{1i}\textbf{u}_{i}+\varrho_{2i}\dfrac{d\textbf{u}_{i}}{dt}+\cdots+\varrho_{s_{i}i}\dfrac{d^{s_{i}-1}\textbf{u}_{i}}{dt^{s_{i}-1}}.
\end{split}
\end{equation}
For $i=1,\cdots,N$, let 
\begin{equation*}
\begin{split}
 \Phi_{i}=\left[
            \begin{array}{cccc}
              0 & 1 & \cdots & 0 \\
              \vdots & \vdots & \ddots & \vdots \\
              0 & 0 & \cdots & 1 \\
              \varrho_{1i} & \varrho_{2i} & \cdots & \varrho_{s_{i}i} \\
            \end{array}
          \right],~ \Gamma_{i}=\left[
                                 \begin{array}{c}
                                   1 \\
                                   0 \\
                                   \vdots \\
                                   0 \\
                                 \end{array}
                               \right]^{T}
\end{split}
\end{equation*}
 and define the following dynamic compensator  \cite{Huang1},
\cite{nik98}:
\begin{equation}\label{doteta1}
\begin{split}
  \dot{\eta}_{i}=M_{i}\eta_{i}+Q_{i}u_{i},
\end{split}
\end{equation}
where $M_{i}\in \mathbb{R}^{s_{i}\times s_{i}}$ is any
Hurwitz matrix and $Q_{i}\in \mathbb{R}^{s_{i}\times 1}$ is any column
vector such that the matrix pair $(M_{i},Q_{i})$ is controllable. As remarked in \cite{Huang1},  \eqref{doteta1} is called a linear internal model of \eqref{system1}, and the following Sylvester equation
\begin{equation}\label{Sylvester}
\begin{split}
  T_{i}\Phi_{i}-M_{i}T_{i}=Q_{i}\Gamma_{i}
\end{split}
\end{equation}
has a unique nonsingular solution $T_{i}$.
For $i=1,\cdots,N$, let
\begin{equation*}
\begin{split}
\Psi_{i}&=\Gamma_{i}T_{i}^{-1}\\
\theta_{i}(v,w)&=T_{i}\mbox{col}(\textbf{u}_{i}(v,w),\dot{\textbf{u}}_{i}(v,w),\cdots,\textbf{u}_{i}^{(s_{i}-1)}(v,w)).\\
\end{split}
\end{equation*}
Then we have
\begin{equation*}
\begin{split}
\frac{\partial \theta_{i}(v,w)}{\partial v}Sv&=(M_{i}+Q_{i}\Psi_{i})\theta_{i}(v,w)\\
\textbf{u}_{i}(v,w)&=\Psi_{i}\theta_{i}(v,w)\\
\end{split}
\end{equation*}
 \end{Remark}
which is called a steady-state generator \cite{Huang1}.
\begin{Remark}
Note that the function $\textbf{u}_{i}(v,w)$  depends on the exogenous signal $v$ and the uncertainty $w$ and thus cannot be directly used to design a feedback control law. As a result, we need to design the linear internal model \eqref{doteta1} to reproduce the function $\theta_{i}(v,w)$ asymptotically and thus provide the information of the function $\textbf{u}_{i}(v,w)$ asymptotically.
\end{Remark}

For $i=1,\cdots,N$ and $k\in\mathbb{Z}^{+}$,  perform the following coordinate and input transformation on  the plant (\ref{system1}) and the internal model
(\ref{doteta1})
\begin{equation}\label{transformation}
\begin{split}
  &\bar{z}_{i}\!=\!z_{i}-\textbf{z}_{i}(v,w),~ \bar{\eta}_{i}=\eta_{i}-\theta_{i}(v,w)-Q_{i}b_{i}^{-1}e_{i}\\
  &e_{i}\!=\!y_{i}\!-\!q(v,w),~ \bar{u}_{i}\!=\!u_{i}\!-\!\Psi_{i}\eta_{i}(t_{k}^{i}),~\forall t\!\in\![t_{k}^{i},t_{k+1}^{i}).
\end{split}
\end{equation}
Then we obtain the  following augmented system
\begin{equation}\label{system2}
\begin{split}
  \dot{\bar{z}}_{i}=&\bar{f}_{i}(\bar{z}_{i},e_{i},\mu)\\
  \dot{\bar{\eta}}_{i}=&M_{i}\bar{\eta}_{i}+M_{i}Q_{i}b_{i}^{-1}e_{i}-Q_{i}b_{i}^{-1}\bar{g}_{i}(\bar{z}_{i},e_{i},\mu)\\
  \dot{e}_{i}=&\bar{g}_{i}(\bar{z}_{i},e_{i},\mu)\!+\!b_{i}\Psi_{i}\bar{\eta}_{i}\!+\!\Psi_{i}Q_{i}e_{i}\!+\!b_{i}\bar{u}_{i}\!+\!b_{i}\Psi_{i}\tilde{\eta}_{i}
\end{split}
\end{equation}
where $i=1,\cdots,N$, $\mu=(v,w)$,
\begin{equation*}
\begin{split}
&\bar{f}_{i}(\bar{z}_{i},e_{i},\mu)=f_{i}(\bar{z}_{i}+\textbf{z}_{i},e_{i}+q,v,w)-f_{i}(\textbf{z}_{i},q,v,w)\\
&\bar{g}_{i}(\bar{z}_{i},e_{i},\mu)=g_{i}(\bar{z}_{i}+\textbf{z}_{i},e_{i}+q,v,w)-g_{i}(\textbf{z}_{i},q,v,w).\\
\end{split}
\end{equation*}
It is easy to verify that,  for any $\mu\in\mathbb{R}^{n_{v}}\times\mathbb{R}^{n_{w}}$ and $i=1,\cdots,N$,
\begin{equation*}
\begin{split}
 \bar{f}_{i}(0,0,\mu)=0,~~\bar{g}_{i}(0,0,\mu)=0.
\end{split}
\end{equation*}
For $i=1,\cdots,N$ and $k\in\mathbb{Z}^{+}$, consider a distributed piecewise constant control law of the following form
  \begin{equation}\label{barui1}
\begin{split}
\bar{u}_{i}(t)=\bar{h}_{i}(e_{vi}(t_{k}^{i})), ~~t\in[t_{k}^{i},t_{k+1}^{i})
\end{split}
\end{equation}
 where $\bar{h}_{i}(\cdot)$ is some globally defined sufficiently smooth function vanishing at the origin.
Denote the state of the closed-loop system composed of the augmented system (\ref{system2}) and the control law (\ref{barui1})  by
$\bar{x}_{c}= \mbox{col } (\bar{z}_{1},\bar{\eta}_{1}, \cdots, \bar{z}_{N},\bar{\eta}_{N},e_{1},\cdots,e_{N})$.

Then we give the following proposition.
\begin{Proposition}\label{Proposition1}
Under Assumptions \ref{Ass2.1}-\ref{Ass2.4}, for any $\epsilon>0$, any known compact sets $\mathbb{V}\in\mathbb{R}^{n_{v}}$ and $\mathbb{W}\in\mathbb{R}^{w}$,  if a distributed control law of the form (\ref{barui1}) can be found such that, for all $\bar{x}_c (0)$, all $v(t)\in\mathbb{V}$ and all $w\in\mathbb{W}$, $\bar{x}_{c}(t)$ exists and is bounded for all $t\geq0$, and satisfies
\begin{equation}\label{pcgs1}
\lim_{t \to \infty}\sup\|\bar{x}_c(t)\|\leq \epsilon,
\end{equation}
then Problem \ref{Problem1} for the system \eqref{system1} is solvable by the following  distributed piecewise constant control law
 \begin{equation}\label{ui2}
\begin{split}
u_{i}(t)&=\bar{h}_{i}({e}_{vi}(t_{k}^{i}))+\Psi_{i}\eta_{i}(t_{k}^{i})\\
\dot{\eta}_{i}(t)&=M_{i}\eta_{i}(t)+Q_{i}u_{i}(t),~~t\in[t_{k}^{i},t_{k+1}^{i})\\
\end{split}
\end{equation}
for $i=1,\cdots,N$ and $k\in\mathbb{Z}^{+}$.
\end{Proposition}
\begin{Proof}
First, 
according to \eqref{pcgs1}, it is easy to see that
\begin{equation}\label{pcgs2}
\lim_{t \to \infty}\sup\|e(t)\|\leq\lim_{t \to \infty}\sup\|\bar{x}_c(t)\|\leq \epsilon.
\end{equation}
Thus Property (2) in Problem \ref{Problem1} is satisfied.

 Next, we only need to show that Property (1) in Problem \ref{Problem1} is also satisfied.
We first denote the state of the closed-loop system composed of the system \eqref{system1} and the control law \eqref{ui2} by $x_{c}=\mbox{col}(z_{1},\eta_{1},\cdots,z_{N},\eta_{N},y_{1},\cdots,y_{N})$. Then, based on the coordinate transformation \eqref{transformation}, we know
 \begin{equation}\label{bxc1}
\begin{split}
x_{c}=&~\bar{x}_{c}+\mbox{col}(\textbf{z}_{1}(v,w),\theta_{1}(v,w)+Q_{1}b_{1}(w)^{-1}e_{1},\\
&\cdots,\textbf{z}_{N}(v,w),\theta_{N}(v,w)+Q_{N}b_{N}(w)^{-1}e_{N},\\
&q(v,w),\cdots,q(v,w)).\\
\end{split}
\end{equation}
Note that, for $i=1,\cdots,N$, $\textbf{z}_{i}(v,w)$, $b_{i}(w)$, $\theta_{i}(v,w)$ and $q(v,w)$ are all smooth functions, and the boundaries of the compact sets $\mathbb{V}$ and $\mathbb{W}$ are known. Then  $\textbf{z}_{i}(v,w)$, $b_{i}(w)$, $\theta_{i}(v,w)$ and $q(v,w)$ are all bounded for all $t\geq0$. 
 Together with \eqref{bxc1} and the fact that $\bar{x}_{c}(t)$ exists and is bounded for all $t\geq0$, we conclude that $x_{c}(t)$ exists and is bounded for all $t\geq0$, i.e., Property (1) in Problem \ref{Problem1} is satisfied.

Thus the proof is completed.
\end{Proof}


We call the problem of designing a control law of the form (\ref{barui1}) to achieve (\ref{pcgs1})  the cooperative global robust practical stabilization problem for the augmented system  \eqref{system2}.

 \begin{Remark}
The transformation \eqref{transformation} is modified from the corresponding one in \cite{DongHuang2014a} by  replacing $\eta_{i}(t)$ with $\eta_{i}(t_{k}^{i})$. This modification is necessary for obtaining a directly implementable digital control law, and it also results in a more complex augmented system \eqref{system2}.
\end{Remark}



\section{Main Result}\label{MR}

In this section, we will consider the cooperative global robust practical stabilization for the augmented system \eqref{system2}. For this purpose,
 we need two more assumptions.
 \begin{Assumption}\label{Ass3.1}
 For $i=1,\cdots,N$, and any compact subset $\Omega \subset \mathbb{R}^{n_{v}}\times\mathbb{R}^{n_{w}}$, there exist some $\mathcal{C}^{1}$ functions $V_{1i}(\bar{z}_{i})$ such that, for any $
 \mu\in\Omega$,  any $\bar{z}_{i}$, and any $e_{i}$,
   \begin{equation}\label{Vzi1}
\begin{split}
\underline{\alpha}_{1i}(\|\bar{z}_{i}\|)\leq
V_{1i}(\bar{z}_{i})\leq\bar{\alpha}_{1i}(\|\bar{z}_{i}\|)
\end{split}
\end{equation}
 \begin{equation}\label{dotVzi1}
\begin{split}
\frac{\partial V_{1i}(\bar{z}_{i})}{\partial \bar{z}_{i}}\bar{f}_{i}(\bar{z}_{i},e_{i},\mu)\leq
-\alpha_{1i}(\|\bar{z}_{i}\|)+\gamma_{1i}(e_{i})
\end{split}
\end{equation}
where $\underline{\alpha}_{1i}(\cdot)$ and $\bar{\alpha}_{1i}(\cdot)$ are  some class $\mathcal{K}_{\infty}$ functions,
 $\alpha_{1i}(\cdot)$ are some known class $\mathcal{K}_{\infty}$ function satisfying $\lim_{s\rightarrow0^{+}}\sup(s^{2}/\alpha_{1i}(s))<\infty$ and  $\gamma_{1i}(\cdot)$ are some known smooth  positive definite functions.
 \end{Assumption}
 \begin{Remark}\label{RemarkAss3.1}
Under Assumption \ref{Ass3.1}, the subsystem $\dot{\bar{z}}_{i}=\bar{f}_{i}(\bar{z}_{i},e_{i},\mu)$ is input-to-state stable (ISS) with $e_{i}$ as the input \cite{Sontag1}. Note that this assumption is a standard assumption for cooperative global robust output regulation problem and can also be found in \cite{DongHuang2014a,SuHuang2015}.
\end{Remark}

\begin{Assumption}\label{Ass3.2}
Every node $i=1,\cdots,N$ is reachable from node $0$ in the digraph $\bar{\mathcal{G}}$.
 \end{Assumption}
 \begin{Remark}\label{RemarkAss3.2}
Assumption \ref{Ass3.2} is weaker than Assumption 3.2 of \cite{DongHuang2014a}, since here we allow the digraph $\bar{\mathcal{G}}$ to be directed. Also, by Lemma 4 of \cite{Hu1}, under Assumption \ref{Ass3.2}, $H$ is an $\mathcal{M}$ matrix. Then, by Theorem 2.5.3 of \cite{Horn1}, there exists a positive definite diagonal matrix $D=\mbox{diag}(d_{1},\cdots,d_{N})$ such that $DH+H^{T}D$ is positive definite.
\end{Remark}

Before giving our main result, we introduce some notation. For    any known compact
subset $\mathbb{W}$ and $i=1,\cdots,N$, there always exist some known positive numbers $b_{m}$ and $b_{M}$ such that,  $b_{m}\leq b_{i}(w)\leq b_{M}$ for all $w\in \mathbb{W}$.
For $i=1,\cdots,N$ and $k\in\mathbb{Z}^{+}$, define
\begin{equation}\label{varthetai1}
\begin{split}
&\vartheta_{i}(t)=-\rho_{i}(e_{vi}(t))e_{vi}(t)\\
&\tilde{\vartheta}_{i}(t)=\vartheta_{i}(t_{k}^{i})-\vartheta_{i}(t),~~\forall t\in[t_{k}^{i},t_{k+1}^{i})\\
\end{split}
\end{equation}
where $\rho_{i}(\cdot)$, $i=1,\cdots,N$, are some sufficiently smooth positive functions to be specified later. Then we consider the following  control law
\begin{equation}\label{barui2}
 \bar{u}_{i}(t)=\vartheta_{i}(t_{k}^{i}),~\forall t\in[t_{k}^{i},t_{k+1}^{i})
\end{equation}
and event-triggered mechanism
\begin{equation}\label{trigger2}
\begin{split}
t_{k+1}^{i}\!\!=\!\inf\{t>t_{k}^{i}~|~ (\tilde{\vartheta}_{i}(t)\!+\!\Psi_{i}\tilde{\eta}_{i}(t))^{2}\!-\!\sigma\vartheta_{i}^{2}(t)\!\geq\!\delta\}
\end{split}
\end{equation}
where $\sigma>0$ and $\delta>0$ are some constants to be determined later. Clearly,  the control law \eqref{barui2} and the event-triggered mechanism \eqref{trigger2} are both distributed and output-based, since the term $\vartheta_{i}(t)$ only depends on $e_{vi}$, i.e., the outputs of the neighbors of agent $i$.
\begin{Remark}
For  any $t\in[t_{k},t_{k+1})$ with $k\in\mathbb{Z}^{+}$, let $\tilde{x}_{c}(t)=\bar{x}_{c}(t_{k})-\bar{x}_{c}(t)$. Then the centralized state-based event-triggered mechanism  in \cite{XieD1} is equivalent to the following form
\begin{equation}\label{trigger0a}
\begin{split}
t_{k+1}\!=\!\inf\{t>t_{k}~|~  \|\tilde{x}_{c}(t)\|\geq\sigma\|\bar{x}_{c}(t)\|\}
\end{split}
\end{equation}
 for some positive constant $\sigma$.
 For $i=1,\cdots,N$, denote   the full state information of the agent $i$ by $x_{i}(t)$ and let $x_{vi}(t)=\sum_{j=1}^{N}h_{ij}x_{i}(t)$. For $t\in[t_{k}^{i},t_{k+1}^{i})$ with $k\in\mathbb{Z}^{+}$, let $\tilde{x}_{vi}(t)=x_{vi}(t_{k}^{i})-x_{vi}(t)$. Then the distributed state-based event-triggered mechanisms in  \cite{Dimarogonas1,Fan1} are equivalent to
 \begin{equation}\label{trigger0b}
\begin{split}
t_{k+1}^{i}\!=\!\inf\{t>t_{k}^{i}~|~  \|\tilde{x}_{vi}(t)\|\geq\sigma\|x_{vi}(t)\|\}
\end{split}
\end{equation}
 for some positive constant $\sigma$,
 the distributed state-based event-triggered mechanisms in  \cite{Cheng1,Seyboth1} are equivalent to
\begin{equation}\label{trigger0c}
\begin{split}
t_{k+1}^{i}\!\!=\!\!\inf\{t\!>\!t_{k}^{i}|\|\tilde{x}_{vi}(t)\|^{2}\!\!\geq\!\sigma\|x_{vi}(t)\|^{2}\!\!+\!\beta \mathbf{e}^{-\alpha t}\!\!+\!\gamma\}
\end{split}
\end{equation}
for some positive constants $\sigma,\alpha,\beta$ and $\gamma$, and  the observer state based event-triggered mechanisms in \cite{HuW2,Zhang1} are equivalent to
\begin{equation}\label{trigger0d}
\begin{split}
t_{k+1}^{i}\!=\!\inf\{t>t_{k}^{i}~|~ \|\varsigma_{i}(t)\|\geq\sigma\|\tilde{\omega}_{i}(t)\|\}
\end{split}
\end{equation}
where $\tilde{\omega}_{i}(t)$ and $\varsigma_{i}(t)$ are the observer state and the observer state measurement error as shown in \cite{Zhang1}, $\sigma$  is some positive constant.

It can be found that the event-triggered mechanisms \eqref{trigger0a} \eqref{trigger0b} and \eqref{trigger0c} are all state-based and the  event-triggered mechanism \eqref{trigger0d} only depends on the state of the observer.  What makes our event-triggered mechanism \eqref{trigger2} different from  \eqref{trigger0a} \eqref{trigger0b}, \eqref{trigger0c} and \eqref{trigger0d} is that \eqref{trigger2} depends not only the output of the plant but also the state of the internal model. As a result, it is more challenging to analyze the stability of the closed-loop system and prevent the Zeno behavior from happening.
\end{Remark}

 According to the event-triggered mechanism \eqref{trigger2}, for $i=1,\cdots,N$ and $k\in\mathbb{Z}^{+}$,  we have
\begin{equation}\label{tildevarthetai1}
\begin{split}
(\tilde{\vartheta}_{i}(t)+\Psi_{i}\tilde{\eta}_{i}(t))^{2}\leq\sigma\vartheta_{i}^{2}(t)+\delta,~\forall t\in[t_{k}^{i},t_{k+1}^{i}).
\end{split}
\end{equation}
According to \eqref{varthetai1}, we know
\begin{equation}\label{varthetai2}
\begin{split}
&\vartheta_{i}(t_{k}^{i})=\tilde{\vartheta}_{i}(t)+\vartheta_{i}(t),~~\forall t\in[t_{k}^{i},t_{k+1}^{i}).\\
\end{split}
\end{equation}
Then, for $i=1,\cdots,N$, the closed-loop system composed of \eqref{system2} and \eqref{barui2} can be written as follows
\begin{equation}\label{system4}
\begin{split}
  &\dot{Z}_{i}=F_{i}(Z_{i},e_{i},\mu)\\
  &\dot{e}_{i}=\tilde{g}_{i}(Z_{i},e_{i},\mu)+b_{i}\vartheta_{i}(t)+b_{i}\tilde{\vartheta}_{i}(t)+b_{i}\Psi_{i}\tilde{\eta}_{i}
\end{split}
\end{equation}
where $Z_{i}=\mbox{col} (\bar{z}_{i},\bar{\eta}_{i})$,
\begin{equation*}
\begin{split}
 F_{i}(Z_{i},e_{i},\mu)=~&\mbox{col} (\bar{f}_{i}(\bar{z}_{i},e_{i},\mu),M_{i}\bar{\eta}_{i}+M_{i}Q_{i}b_{i}^{-1}e_{i}\\
 &\!-\!Q_{i}b_{i}^{-1}\bar{g}_{i}(\bar{z}_{i},e_{i},\mu))\\
  \tilde{g}_{i}(Z_{i},e_{i},\mu) =~& \bar{g}_{i}(\bar{z}_{i},e_{i},\mu)+b_{i}\Psi_{i}\bar{\eta}_{i}+\Psi_{i}Q_{i}e_{i}.\\
\end{split}
\end{equation*}
Let $Z=\mbox{col}(Z_{1},\cdots,Z_{N})$, $\tilde{\eta}=\mbox{col}(\tilde{\eta}_{1},\cdots,\tilde{\eta}_{N})$, $B=$ $\mbox{diag}(b_{1},\cdots,b_{N})$, $\Psi=\mbox{diag}(\Psi_{1},\cdots,\Psi_{N}),
\vartheta=\mbox{col}(\vartheta_{1},\cdots$, $\vartheta_{N})$, and $\tilde{\vartheta}=\mbox{col}(\tilde{\vartheta}_{1},\cdots,\tilde{\vartheta}_{N})$. 
Then the closed-loop system \eqref{system4} can be put into the following form
\begin{equation}\label{system5}
\begin{split}
  &\dot{Z}=F(Z,e,\mu)\\
  &\dot{e}=\tilde{g}(Z,e,\mu)+B\vartheta+B\tilde{\vartheta}+B\Psi\tilde{\eta}\\
\end{split}
\end{equation}
where
\begin{equation*}
\begin{split}
&F(Z,e,\mu)=\mbox{col}(F_{1}(Z_{1},e_{1},\mu),\cdots,F_{N}(Z_{N},e_{N},\mu))\\
&\tilde{g}(Z,e,\mu)=\mbox{col}(\tilde{g}_{1}(Z_{1},e_{1},\mu),\cdots,\tilde{g}_{N}(Z_{i},e_{N},\mu)).
\end{split}
\end{equation*}
For convenience, we further put \eqref{system5} into the following compact form:
\begin{equation}\label{system6}
\begin{split}
  &\dot{\bar{x}}_{c}=f_{c}(\bar{x}_{c},\mu)\\
\end{split}
\end{equation}
where $\bar{x}_{c}$ is the same as that defined before and $f_{c}(\bar{x}_{c},\mu)=\mbox{col}(F(Z,e,\mu),\tilde{g}(Z,e,\mu)+B\vartheta+B\tilde{\vartheta}+B\Psi\tilde{\eta})$.

Suppose that the solution $\bar{x}_{c}(t)$ of \eqref{system6} under the event-triggered mechanism \eqref{trigger2} is right maximally defined for all $t\in[0,T_{M})$ with $0<T_{M}\leq\infty$.
Then we give the following lemma.
 \begin{Lemma}\label{Lemma1}
 Under Assumptions \ref{Ass2.1}-\ref{Ass3.2},  for $i=1,\cdots,N$, let $\rho_{i}(e_{vi})=a\omega_{i}(e^{2}_{vi})$ where $a$ is a positive real number and $\omega_{i}(\cdot)\geq1$ are some smooth 
  functions. Then, there exists a $\mathcal{C}^{1}$ function $U(\bar{x}_{c})$ and two class $\mathcal{K}_{\infty}$ functions $\underline{\beta}(\cdot)$ and $\bar{\beta}(\cdot)$, such that, for any $\mu\in\Omega$ and any $\bar{x}_{c}$, 
 \begin{equation}\label{U1}
\begin{split}
\underline{\beta}(\|\bar{x}_{c}\|)\leq
U(\bar{x}_{c})\leq\bar{\beta}(\|\bar{x}_{c}\|)
\end{split}
\end{equation}
\begin{equation}\label{dotU4}
\begin{split}
\frac{\partial U(\bar{x}_{c})}{\partial \bar{x}_{c}}f_{c}(\bar{x}_{c},\mu)\leq -\|\bar{x}_{c}\|^{2},~~ \forall~\|\bar{x}_{c}\|\geq\sqrt{\lambda_{3}N\delta}\\
\end{split}
\end{equation}
where $\lambda_{3}=b_{M}^{2}\|DH\|$.
\end{Lemma}


\begin{Proof}
 First, under Assumption \ref{Ass3.1}, by applying Lemma 3.1 of \cite{XuHuang2010}, for $i=1,\cdots,N$, there exist some $\mathcal{C}^{1}$ functions $V_{2i}(Z_{i})$ such that, for any $\mu\in\Omega$, any $Z_{i}$, and any $e_{i}$,
   \begin{equation}\label{VZi2}
\begin{split}
\underline{\alpha}_{2i}(\|Z_{i}\|)\leq
V_{2i}(Z_{i})\leq\bar{\alpha}_{2i}(\|Z_{i}\|)
\end{split}
\end{equation}
 \begin{equation}\label{dotVZi2}
\begin{split}
\frac{\partial V_{2i}(Z_{i})}{\partial Z_{i}}F_{i}(Z_{i},e_{i},\mu)\leq-\|Z_{i}\|^{2}+\gamma_{2i}(e_{i})
\end{split}
\end{equation}
where $\underline{\alpha}_{2i}(\cdot)$ and $\bar{\alpha}_{2i}(\cdot)$ are some class $\mathcal{K}_{\infty}$ functions,  and  $\gamma_{2i}(\cdot)$ are some known smooth  positive definite functions.

Then, by applying the changing supply pair technique \cite{Sontag1}, for  $i=1,\cdots,N$ and any given smooth function $\Delta_{i}(Z_{i})>0$, there exist some $\mathcal{C}^{1}$ functions $V_{3i}(Z_{i})$ such that, for any $\mu\in\Omega$, any $Z_{i}$, and any $e_{i}$, 
   \begin{equation}\label{VZi3}
\begin{split}
\underline{\alpha}_{3i}(\|Z_{i}\|)\leq
V_{3i}(Z_{i})\leq\bar{\alpha}_{3i}(\|Z_{i}\|)
\end{split}
\end{equation}
 \begin{equation}\label{dotVZi3}
\begin{split}
\frac{\partial V_{3i}(Z_{i})}{\partial Z_{i}}F_{i}(Z_{i},e_{i},\mu)\leq
-\Delta_{i}(Z_{i})\|Z_{i}\|^{2}+\pi_{i}(e_{i})e_{i}^{2}
\end{split}
\end{equation}
where $\underline{\alpha}_{3i}(\cdot)$ and $\bar{\alpha}_{3i}(\cdot)$ are some class $\mathcal{K}_{\infty}$ functions,  and  $\pi_{i}(\cdot)$ are some known smooth  positive functions. Let $V_{3}(Z)=\sum_{i=1}^{N}V_{3i}(Z_{i})$. Then, for any $\mu\in\Omega$, any $Z$, and any $e$, we have  
   \begin{equation}\label{VZ3}
\begin{split}
\underline{\alpha}_{3}(\|Z\|)\leq
V_{3}(Z)\leq\bar{\alpha}_{3}(\|Z\|)
\end{split}
\end{equation}
 \begin{equation}\label{dotVZ3}
\begin{split}
\frac{\partial V_{3}(Z)}{\partial Z}F(Z,e,\mu)\leq& \sum_{i=1}^{N}(-\Delta_{i}(Z_{i})\|Z_{i}\|^{2}+ \pi_{i}(e_{i})e_{i}^{2})
\end{split}
\end{equation}
for some class $\mathcal{K}_{\infty}$ functions $\underline{\alpha}_{3}(\cdot)$ and $\bar{\alpha}_{3}(\cdot)$.

According to \eqref{system5}, we have
\begin{equation}\label{dotev1}
\begin{split}
\dot{e}_{v}&=H\dot{e}\\
&=H\tilde{g}(Z,e,\mu)+HB(\vartheta(t)+\tilde{\vartheta}(t)+\Psi\tilde{\eta}(t))\\
&=G(Z,e_{v},\mu)+HB(\vartheta(t)+\tilde{\vartheta}(t)+\Psi\tilde{\eta}(t))\\
\end{split}
\end{equation}
where $G(Z,e_{v},\mu)=\mbox{col}(G_{1}(Z,e_{v},\mu),\cdots,G_{N}(Z,e_{v},\mu))=H\tilde{g}(Z,e,\mu)$. Let $H=\mbox{col}(H_{1},\cdots,H_{N})$ where $H_{i}=[h_{i1}, \cdots,h_{iN}]$ for $i=1,\cdots,N$.
As noted in Remark \ref{RemarkAss3.2}, 
$DH+H^{T}D$ is a positive definite matrix, thus $B(DH$ $+H^{T}D)B$ is also a positive definite matrix. Define $\lambda_{1}=b_{m}^{2}\lambda_{\min}(DH+H^{T}D)$,
$\lambda_{2}=b_{M}\|DH\|$. 

Next, motivated from equation (37) of \cite{SuHuang2015}, let
\begin{equation}\label{V4e1}
\begin{split}
V_{4}(e)=\sum_{i=1}^{N}d_{i}b_{i}\int_{0}^{e^{2}_{vi}}\omega_{i}(s)ds.
\end{split}
\end{equation}
It can be seen that $V_{4}(e)$ is positive definite and radially unbounded. Thus,  
there  exist some class $\mathcal{K}_{\infty}$ functions $\underline{\alpha}_{4}(\cdot)$ and $\bar{\alpha}_{4}(\cdot)$ such that
\begin{equation}\label{V4e2}
\begin{split}
\underline{\alpha}_{4}(\|e\|)\leq V_{4}(e) \leq\bar{\alpha}_{4}(\|e\|).
\end{split}
\end{equation}
Let $e_{v}^{*}=\mbox{col}(\omega_{1}(e^{2}_{v1})e_{v1},\cdots,\omega_{N}(e^{2}_{vN})e_{vN})$.
Then, by \eqref{dotev1} and \eqref{V4e1},  for any $\mu\in\Omega$, and any $e_{v}$,  
\begin{equation}\label{dotV4e1}
\begin{split}
&\frac{\partial V_{4}(e)}{\partial e}(\tilde{g}(Z,e,\mu)+B\vartheta+B\tilde{\vartheta}+B\Psi\tilde{\eta})\\
=&\frac{\partial V_{4}(H^{-1}e_{v})}{\partial e_{v}}\!( G(Z,e_{v},\mu)\!\!+\!\!H\!B(\vartheta(t)\!+\!\tilde{\vartheta}(t)\!+\!\Psi\tilde{\eta}(t)))\\
=&2\sum_{i=1}^{N}d_{i}b_{i}\omega_{i}(e^{2}_{vi})e_{vi}\dot{e}_{vi}\\
 =&2\sum_{i=1}^{N}d_{i}b_{i}\omega_{i}(e^{2}_{vi})e_{vi}\big(G_{i}(Z,e_{v},\mu)\\
 &+H_{i}B (\vartheta(t)+\tilde{\vartheta}(t)+\Psi\tilde{\eta})\big)\\
 =& 2(e_{v}^{*})^{T}BD\big(G(Z,e_{v},\mu)\!+\!H\!B (\vartheta(t)\!+\!\tilde{\vartheta}(t)\\
 &+\Psi\tilde{\eta}(t))\big).
\end{split}
\end{equation}
Note that 
\begin{equation}\label{inequality1}
\begin{split}
&2(e_{v}^{*})^{T}BDG(Z,e_{v},\mu)\\
=&2(e_{v}^{*})^{T}BDH\tilde{g}(Z,e,\mu)\\
\leq&2\|e_{v}^{*}\|\|B\|\|DH\|\|\tilde{g}(Z,e,\mu)\|\\
\leq&2\lambda_{2}\|e_{v}^{*}\|\|\tilde{g}(Z,e,\mu)\|\\
\leq &\lambda_{2}(\|e_{v}^{*}\|^{2}+\|\tilde{g}(Z,e,\mu)\|^{2})\\
\end{split}
\end{equation}
\begin{equation}\label{inequality2}
\begin{split}
&2(e_{v}^{*})^{T}BDHB \vartheta(t)\\
=&-2a(e_{v}^{*})^{T}BDHBe_{v}^{*}\\
=&-a(e_{v}^{*})^{T}B(DH+H^{T}D)Be_{v}^{*}\\
\leq&-a\lambda_{1}\|e_{v}^{*}\|^{2}.
\end{split}
\end{equation}
From \eqref{varthetai1} and \eqref{tildevarthetai1},  we know
\begin{equation}\label{inequality3}
\begin{split}
&2(e_{v}^{*})^{T}BDHB(\tilde{\vartheta}(t)+\Psi\tilde{\eta}(t))\\
\leq&2\|e_{v}^{*}\|\|B\|\|DH\|\|B\|\|\tilde{\vartheta}(t)+\Psi\tilde{\eta}(t)\|\\
\leq&2\lambda_{3}\|e_{v}^{*}\|\|\tilde{\vartheta}(t)+\Psi\tilde{\eta}(t)\|\\
\leq&\lambda_{3}(\|e_{v}^{*}\|^{2}+\|\tilde{\vartheta}(t)+\Psi\tilde{\eta}(t)\|^{2})\\
\leq&\lambda_{3}(\|e_{v}^{*}\|^{2}+\sigma\|\vartheta(t)\|^{2}+N\delta)\\
=&\lambda_{3}(\|e_{v}^{*}\|^{2}+\sigma a^2 \|e_{v}^{*}\|^{2}+N\delta)\\
=&\lambda_{3}(1+\sigma a^2) \|e_{v}^{*}\|^{2}+ \lambda_{3} N\delta.
\end{split}
\end{equation}
Also, since $\tilde{g}_{i}(Z_{i},e_{i},\mu)$
is smooth and $\tilde{g}_{i}(0,0,\mu)=0$ for all $\mu\in\Omega$, by
Lemma 7.8 of \cite{Huang1}, there exist some smooth functions $\bar{\varphi}_{i}(Z_{i})$ and $\bar{\chi}_{i}(e_{i})$ satisfying $\bar{\varphi}_{i}(0)=0$ and $\bar{\chi}_{i}(0)=0$, such
that, for any $Z_{i}\in\mathbb{R}^{n+s_{i}}$, $e_{i}\in \mathbb{R}$, and $\mu\in\Omega$,
\begin{equation}\label{tildegi1}
\begin{split}
  |\tilde{g}_{i}(Z_{i},e_{i},\mu)|\leq&\bar{\varphi}_{i}(Z_{i})+\bar{\chi}_{i}(e_{i})
\end{split}
\end{equation}
which implies
\begin{equation}\label{tildegi2}
\begin{split}
  &|\tilde{g}_{i}(Z_{i},e_{i},\mu)|^{2}\\
  \leq&(\bar{\varphi}_{i}(Z_{i})+\bar{\chi}_{i}(e_{i}))^{2}\\
  =&|\bar{\varphi}_{i}(Z_{i})|^{2}+2\bar{\varphi}_{i}(Z_{i})\bar{\chi}_{i}(e_{i})+|\bar{\chi}_{i}(e_{i})|^{2}\\
  \leq&2|\bar{\varphi}_{i}(Z_{i})|^{2}+2|\bar{\chi}_{i}(e_{i})|^{2}.\\
\end{split}
\end{equation}
Since $\bar{\varphi}_{i}(0)=0$ and $\bar{\chi}_{i}(0)=0$, there exist two smooth positive functions $\varphi_{i}(\cdot)$ and $\chi_{i}(\cdot)$ such that, for any $Z_{i}\in\mathbb{R}^{n+s_{i}}$, $e_{i}\in \mathbb{R}$,
\begin{equation}\label{bvarphi1}
\begin{split}
2|\bar{\varphi}_{i}(Z_{i})|^{2}\leq&\varphi_{i}(Z_{i})\|Z_{i}\|^{2}\\
2|\bar{\chi}_{i}(e_{i})|^{2}\leq&\chi_{i}(e_{i})e_{i}^{2}.\\
\end{split}
\end{equation}
Combining \eqref{tildegi2} and \eqref{bvarphi1}, we have
\begin{equation}\label{tildegi3}
\begin{split}
  |\tilde{g}_{i}(Z_{i},e_{i},\mu)|^{2} \leq&\varphi_{i}(Z_{i})\|Z_{i}\|^{2}+\chi_{i}(e_{i})e_{i}^{2}.\\
\end{split}
\end{equation}
 Then, according to  \eqref{dotV4e1}, \eqref{inequality1}, \eqref{inequality2}, \eqref{inequality3} and \eqref{tildegi3},  for any $\mu\in\Omega$, any $Z_{i}$ and $e_{i}$,  we have
\begin{equation}\label{dotV4e2}
\begin{split}
&\frac{\partial V_{4}(e)}{\partial e}(\tilde{g}(Z,e,\mu)+B\vartheta+B\tilde{\vartheta}+B\Psi\tilde{\eta})\\
\leq&~\lambda_{2}(\|e_{v}^{*}\|^{2}+\|\tilde{g}(Z,e,\mu)\|^{2})-a\lambda_{1}\|e_{v}^{*}\|^{2}\\
&+\lambda_{3}(1+\sigma a^2) \|e_{v}^{*}\|^{2}+ \lambda_{3} N\delta \\
\leq& -\sum_{i=1}^{N}(a\lambda_{1}-\lambda_{2}-\lambda_{3}-\lambda_{3}\sigma a^{2})\omega_{i}^{2}(e^{2}_{vi})e^{2}_{vi}\\
&+\sum_{i=1}^{N}\lambda_{2}(\varphi_{i}(Z_{i})\|Z_{i}\|^{2}+\chi_{i}(e_{i})e_{i}^{2})+\lambda_{3}N\delta.\\
\end{split}
\end{equation}
Let $U(\bar{x}_{c})=V_{3}(Z)+V_{4}(e)$. 
Clearly, there exist two class $\mathcal{K}_{\infty}$ functions $\underline{\beta}(\cdot)$ and $\bar{\beta}(\cdot)$ such that
(\ref{U1}) is satisfied.

Also, according to \eqref{dotVZ3} and \eqref{dotV4e2}, for any $\mu\in\Omega$, and any $\bar{x}_{c}$,  we have
\begin{equation}\label{dotU2}
\begin{split}
&\frac{\partial U(\bar{x}_{c}) }{\partial \bar{x}_{c}} f_{c}(\bar{x}_{c},\mu)\\
\leq &-\sum_{i=1}^{N}(\Delta_{i}(Z_{i})-\lambda_{2}\varphi_{i}(Z_{i}))\|Z_{i}\|^{2}\\
&-\sum_{i=1}^{N}(a\lambda_{1}-\lambda_{2}-\lambda_{3}-\lambda_{3}\sigma a^{2})\omega_{i}^{2}(e^{2}_{vi})e^{2}_{vi}\\
&+\hat{\rho}(e)+\lambda_{3}N\delta\\
\end{split}
\end{equation}
with $\hat{\rho}(e)=\sum_{i=1}^{N}\big(\lambda_{2}\chi_{i}(e_{i})+\pi_{i}(e_{i})\big)e_{i}^{2}$.
Let
\begin{equation}\label{barrho1}
\begin{split}
\bar{\rho}(e_{v})=\sqrt{\hat{\rho}(H^{-1}e_{v})}=\sqrt{\hat{\rho}(e)}.
\end{split}
\end{equation}
Note that $\bar{\rho}(0)=0$. Then, by applying lemma $7.8$ of \cite{Huang1} again,
there exist some smooth  functions $\check{\rho}_{i}(e_{vi})$ with $\check{\rho}_{i}(0)=0$ and  $i=1,\cdots,N$, such that, for any $e_{vi}\in\mathbb{R}$,
 \begin{equation}\label{barrho2}
\begin{split}
\bar{\rho}(e_{v})\leq\sum_{i=1}^{N}\check{\rho}_{i}(e_{vi})
\end{split}
\end{equation}
which in turn implies
 \begin{equation}\label{hatrho1}
\begin{split}
\hat{\rho}(e)=|\bar{\rho}(e_{v})|^{2}\leq&|\sum_{i=1}^{N}\check{\rho}_{i}(e_{vi})|^{2}\\
\leq&\sum_{i=1}^{N}N|\check{\rho}_{i}(e_{vi})|^{2}.\\
\end{split}
\end{equation}
Since $\check{\rho}_{i}(0)=0$, there exist some smooth positive functions $\tilde{\rho}_{i}(e_{vi})\geq1$ with $i=1,\cdots,N$, such that, for any $e_{vi}\in\mathbb{R}$,
\begin{equation}\label{checkrho1}
\begin{split}
N|\check{\rho}_{i}(e_{vi})|^{2}\leq\tilde{\rho}_{i}(e_{vi})e_{vi}^{2}.
\end{split}
\end{equation}
According to \eqref{hatrho1} and \eqref{checkrho1}, we have
 \begin{equation}\label{hatrho2}
\begin{split}
\hat{\rho}(e)\leq&\sum_{i=1}^{N}\tilde{\rho}_{i}(e_{vi})e_{vi}^{2}.\\
\end{split}
\end{equation}
Choose
 \begin{equation}\label{Delta1}
\begin{split}
&\Delta_{i}(Z_{i})\geq\lambda_{2}\varphi_{i}(Z_{i})+2\\
&a\geq \dfrac{1}{\lambda_{1}}(\lambda_{2}+2\lambda_{3}+1)\\
&0<\sigma\leq\frac{1}{a^{2}}.\\
\end{split}
\end{equation}
Then, together with \eqref{dotU2} and \eqref{hatrho2}, we have
\begin{equation}\label{dotU3}
\begin{split}
&\frac{\partial U(\bar{x}_{c}) }{\partial \bar{x}_{c}} f_{c}(\bar{x}_{c},\mu)\\
\leq &-2\sum_{i=1}^{N}\|Z_{i}\|^{2}-\sum_{i=1}^{N}(\omega_{i}^{2}(e^{2}_{vi})-\tilde{\rho}_{i}(e_{vi}))e^{2}_{vi}\\
&+\lambda_{3}N\delta.\\
\end{split}
\end{equation}
Let $\omega_{i}(\cdot)$ be a smooth function satisfying
\begin{equation}\label{omegai1}
\begin{split}
\omega_{i}^{2}(e_{vi}^{2})\geq\omega_{i}(e_{vi}^{2})\geq \tilde{\rho}_{i}(e_{vi})+\frac{2}{\lambda_{\min}(H^{2})}.
\end{split}
\end{equation}
Then we have
\begin{equation}\label{dotU5}
\begin{split}
&\frac{\partial U(\bar{x}_{c}) }{\partial \bar{x}_{c}} f_{c}(\bar{x}_{c},\mu)\\
\leq &-2\sum_{i=1}^{N}\|Z_{i}\|^{2}-\sum_{i=1}^{N}\frac{2e^{2}_{vi}}{\lambda_{\min}(H^{2})}+\lambda_{3}N\delta\\
\leq&-2\|Z\|^{2}-2\|e\|^{2}+\lambda_{3}N\delta\\
=&-2\|\bar{x}_{c}\|^{2}+\lambda_{3}N\delta\\
\leq& -\|\bar{x}_{c}\|^{2},~~\forall~\|\bar{x}_{c}\|\geq\sqrt{\lambda_{3}N\delta}.
\end{split}
\end{equation}
Thus the proof is completed.
\end{Proof}
\begin{Remark}\label{RemarkLemma1}
Lemma \ref{Lemma1}  implies that, for any  $t\in[0,T_M)$,
\begin{equation*}
\begin{split}
\|\bar{x}_{c} (t)\| \leq \max \{\sqrt{\lambda_{3}N\delta}, \underline{\beta}^{-1}(\bar{\beta}(\|\bar{x}_{c}(0)\|))\}.
\end{split}
\end{equation*}
That is to say, $\bar{x}_{c}(t)$ is bounded over $[0,T_M)$.
\end{Remark}

 Lemma \ref{Lemma1} together with Proposition \ref{Proposition1} leads to our main result as follows.

\begin{Theorem}\label{Theorem1}
 Under Assumptions \ref{Ass2.1}-\ref{Ass3.2}, for $k\in\mathbb{Z}^{+}$, $i=1,\cdots,N$,  and any $\epsilon > 0$, 
the cooperative global robust practical output regulation problem for the system \eqref{system1} is solvable by the following distributed output feedback control law
\begin{equation}\label{ui3}
\begin{split}
 u_{i}(t)&=-\rho_{i}(e_{vi}(t_{k}^{i}))e_{vi}(t_{k}^{i})+\Psi_{i}\eta_{i}(t_{k}^{i}) \\
 \dot{\eta}_{i}(t)&=M_{i}\eta_{i}(t)+Q_{i}u_{i}(t),~\forall t\in[t_{k}^{i},t_{k+1}^{i})
 \end{split}
\end{equation}
under the distributed output-based event-triggered mechanism \eqref{trigger2} with  $\delta = \frac{(\bar{\beta}^{-1} (\underline{\beta}( \epsilon)))^2}{\lambda_{3}N}$.
\end{Theorem}
\begin{Proof}
We first consider the case that, for all $i=1,\cdots,N$, the number of the triggering times is finite. Then, there exists a finite time $T_{0}$ such that the closed-loop system \eqref{system6} is a nonlinear time-invariant continuous-time system for all $t\geq T_{0}$. Thus, together with Remark \ref{RemarkLemma1}, we conclude that $T_{M}=\infty$. 

Next, we consider the case that, for any $i=1,\cdots,N$,  the time sequence $\{t_{k}^{i}\}$ has  infinite members. In this case, if we show $\lim_{k\rightarrow\infty}t_{k}^{i}=\infty$, then $T_{M}$ must be equal to $\infty$.
Note that, for any $t\in[0, T_{M})$, 
 \begin{equation}\label{dottildevarthetai1}
\begin{split}
&\frac{d(\tilde{\vartheta}_{i}(t)+\Psi_{i}\tilde{\eta}_{i}(t))^{2}}{dt}\\
=&2(\tilde{\vartheta}_{i}(t)+\Psi_{i}\tilde{\eta}_{i}(t))(\dot{\tilde{\vartheta}}_{i}(t)+\Psi_{i}\dot{\tilde{\eta}}_{i}(t))\\
=&2(\tilde{\vartheta}_{i}(t)+\Psi_{i}\tilde{\eta}_{i}(t))(-\dot{\vartheta}_{i}(t)-\Psi_{i}\dot{\eta}_{i}(t))\\
=&-2(\tilde{\vartheta}_{i}(t)+\Psi_{i}\tilde{\eta}_{i}(t))\bigg(\big(\frac{d\rho_{i}(e_{vi})}{de_{vi}}e_{vi}\\
&+\rho_{i}(e_{vi})\big)\dot{e}_{vi}(t)+\Psi_{i}\dot{\eta}_{i}(t)\bigg).\\
\end{split}
\end{equation}
Also, for any $t\in[0, T_{M})$, we know
 \begin{equation}\label{evi2}
\begin{split}
&e_{vi}(t)=\sum_{j=1}^{N}h_{ij}e_{j}(t)\\
&\tilde{\vartheta}_{i}(t)=\rho_{i}(e_{vi}(t))e_{vi}(t)-\rho_{i}(e_{vi}(t_{k}^{i}))e_{vi}(t_{k}^{i})\\
&\tilde{\eta}_{i}(t)=\eta_{i}(t_{k}^{i})-\eta_{i}(t)\\
&\dot{e}_{vi}(t)\!=\sum_{j=1}^{N}h_{ij}\dot{e}_{j}(t)\\
&~~~~~~~=\!\!\sum_{j=1}^{N}h_{ij}(\tilde{g}_{j}(Z_{j},e_{j},\mu)\!-\!b_{j}\rho_{j}(e_{vj}(t_{k}^{j}))e_{vj}(t_{k}^{j})\\
&~~~~~~~~~~+b_{j}\Psi_{j}\tilde{\eta}_{j}(t))\\
&\dot{\eta}_{i}(t)\!=\!M_{i}\eta_{i}(t)\!+\!Q_{i}(\Psi_{i}\eta_{i}(t_{k}^{i})\!-\!\rho_{i}(e_{vi}(t_{k}^{i}))e_{vi}(t_{k}^{i})).\\
\end{split}
\end{equation}
Since, by \eqref{bxc1} and Remark \ref{RemarkLemma1}, the state $x_{c}(t)$ and thus the states $z_{i}(t)$, $\eta_{i}(t)$ and the tracking error $e_{i}(t)$ of the closed-loop system composed of \eqref{system1} and \eqref{ui3} are bounded for all  $t\in[0,T_{M})$.
Thus, \eqref{evi2} implies  that $\tilde{\vartheta}_{i}(t)$, $\tilde{\eta}_{i}(t)$, $\dot{e}_{vi}(t)$ and $\dot{\eta}_{i}(t)$ are all bounded for all $t\in[0,T_{M})$.
Thus, 
there always exists a positive number $c_{0}$ depending on $\delta$ and $\bar{x}_{c}(0)$, 
 such that 
 \begin{equation}\label{dottildevarthetai2}
\begin{split}
\bigg|\frac{d(\tilde{\vartheta}_{i}(t)+\Psi_{i}\tilde{\eta}_{i}(t))^{2}}{dt}\bigg|\leq c_{0},~\forall t\in[0,T_{M}).
\end{split}
\end{equation}

On the other hand, from the second equation of \eqref{tildeevi1} and the second equation of \eqref{varthetai1}, we know
 \begin{equation}\label{tildevarthetai2}
\begin{split}
&\tilde{\eta}_{i}(t_{k}^{i})=\eta_{i}(t_{k}^{i})-\eta_{i}(t_{k}^{i})=0\\
&\tilde{\vartheta}_{i}(t_{k}^{i})=\vartheta_{i}(t_{k}^{i})-\vartheta_{i}(t_{k}^{i}=0,
\end{split}
\end{equation}
and, from \eqref{trigger2}, we have
 \begin{equation}\label{tildevarthetai3}
\begin{split}
&\lim_{t\rightarrow (t_{k+1}^{i})^{-}}\!\!(\tilde{\vartheta}_{i}(t)+\Psi_{i}\tilde{\eta}_{i}(t))^{2}\geq\lim_{t\rightarrow (t_{k+1}^{i})^{-}}\!\!(\sigma\vartheta_{i}^{2}(t)+\delta)\geq\delta.\\
\end{split}
\end{equation}
Combining \eqref{dottildevarthetai2}, \eqref{tildevarthetai2} and \eqref{tildevarthetai3}, we can conclude that, for any $i = 1, \cdots, N$,  $k \in \mathbb{Z}^+$, 
$t_{k+1}^{i} - t_{k}^{i} \geq \frac{\delta}{c_{0}}$, that is to say, $\tau_{d}\geq \frac{\delta}{c_{0}}$.
Together with the condition that the sequence $\{t_{k}^{i}\}$ has infinite members, we have $\lim_{k\rightarrow\infty}t_{k}^{i}=\infty$, which further implies that the  event-triggered mechanism \eqref{trigger2} does not exhibit the Zeno behavior. As a result, the solution $\bar{x}_{c}(t)$ of the closed-loop system \eqref{system6} must exist for all time, i.e., $T_{M}=\infty$.

Since the solution $\bar{x}_{c}(t)$ of  \eqref{system6} exists  for all $t\in[0,\infty)$ whether the number of the triggering numbers is finite or infinite,  by Theorem 4.18 of \cite{Khalil1} and Lemma \ref{Lemma1} here, we have that  the solution of the closed-loop system \eqref{system6} is globally ultimately bounded with the ultimate bound $d(\delta)=\underline{\beta}^{-1}(\bar{\beta}(\sqrt{\lambda_{3}N\delta}))$, that is,
 \begin{equation}\label{barxc1}
\begin{split}
\lim_{t\rightarrow\infty}\sup\|\bar{x}_{c}(t)\|\leq d(\delta)= \underline{\beta}^{-1}(\bar{\beta}(\sqrt{\lambda_{3}N\delta})).\\
\end{split}
\end{equation}
Note that $d(\cdot)$ is  an invertible class $\mathcal{K}_{\infty}$ function, since both $\underline{\beta}^{-1}(\cdot)$ and $\bar{\beta}(\cdot)$ are invertible class $\mathcal{K}_{\infty}$ functions.
Then, for any $\epsilon > 0$, letting
 \begin{equation}\label{delta1}
\begin{split}
 \delta=d^{-1}(\epsilon)=\frac{(\bar{\beta}^{-1} (\underline{\beta}( \epsilon)))^2}{\lambda_{3}N}\\
\end{split}
\end{equation}
  gives $\lim_{t\rightarrow\infty}\sup\|\bar{x}_{c}(t)\|\leq\epsilon$.
 By applying Proposition \ref{Proposition1}, the proof is thus completed.
%
%
\end{Proof}

\begin{Remark}
It is of interest to discuss a special case of Theorem \ref{Theorem1} by letting $\delta=0$. In this case, according to \eqref{barxc1}, we can conclude that
\begin{equation}\label{et1}
\begin{split}
\lim_{t\rightarrow\infty}\|e(t)\| = \lim_{t\rightarrow\infty}\|\bar{x}_{c}(t)\|=0\\
\end{split}
\end{equation}
which means that the tracking error approaches zero asymptotically as time tends to infinity.
However, under such a case, the existence of the minimal inter-execution time cannot be guaranteed  from the inequality $\tau_{d}\geq \frac{\delta}{c_{0}}$, and thus the event-triggered mechanism may exhibit the Zeno behavior.    Due to this reason, we have introduced the parameter $\delta$ in the event-triggered mechanism \eqref{trigger2}.
What is more, as will be observed from the simulation example in next section, a larger $\delta$ usually  leads to larger steady tracking error but leads to less triggering number.
\end{Remark}
\begin{Remark}\label{Discretization}
Note that the control law \eqref{ui3} directly leads to the following digital implementation:
\begin{equation}\label{ui4}
\begin{split}
 u_{i}(t)=&-\rho_{i}(e_{vi}(t_{k}^{i}))e_{vi}(t_{k}^{i})+\Psi_{i}\eta_{i}(t_{k}^{i}) \\
 \eta_{i}(t_{k+1}^{i})\!\!=&\mathbf{e}^{M_{i}(t_{k+1}^{i}-t_{k}^{i})}\eta_{i}(t_{k}^{i})+Q_{i}\big(\Psi_{i}\eta_{i}(t_{k}^{i})-\\
 &\rho_{i}(e_{vi}(t_{k}^{i}))e_{vi}(t_{k}^{i})\big)\int_{t_{k}^{i}}^{t_{k+1}^{i}}\mathbf{e}^{M_{i}(t_{k+1}^{i}-\tau)}d\tau.\\
 \end{split}
\end{equation}

In contrast,  in some existing literature on the event-triggered cooperative control problems such as \cite{Wang5}  or \cite{Zhang1},  the control laws are not piecewise constant and thus cannot be directly implemented in a digital platform. For example,    the control law in \cite{Zhang1} takes the following form:   
\begin{equation}\label{Au2}
\begin{split}
u_{i}(t)&=G_{1}\phi_{i}(t),~ t\in[t_{k}^{i},t_{k+1}^{i})\\
\dot{\phi}_{i}(t)&=G_{2}\phi_{i}(t)+G_{3}\phi_{vi}(t_{k}^{i})+G_{4}e_{vi}(t_{k}^{i})\\
\end{split}
\end{equation}
where $k\in\mathbb{Z}^{+}$, $G_{1},G_{2},G_{3}$ and $G_{4}$ are some matrices with proper dimensions, $\phi_{i}$ subsystem is a dynamic compensator and $\phi_{vi}(t)=\sum_{j=1}^{N}h_{ij}\phi_{i}(t)$.
It can be seen that the signal generated by the control law \eqref{Au2} is not piecewise constant. To implement
\eqref{Au2}  in a digital platform, one has to further sample \eqref{Au2} to obtain the following:
\begin{equation}\label{Au3}
\begin{split}
&u_{i}(t)=G_{1}\phi_{i}(t_{k}^{i}), ~ t\in[t_{k}^{i},t_{k+1}^{i}) \\
&{\phi}_{i}(t_{k+1}^{i})=\mathbf{e}^{G_{2}(t_{k+1}^{i}-t_{k}^{i})}\phi_{i}(t_{k}^{i})+\big(G_{3}\phi_{vi}(t_{k}^{i})\\
&~~~~~~~~~~~~~~+G_{4}e_{vi}(t_{k}^{i})\big)\int_{t_{k}^{i}}^{t_{k+1}^{i}}\!\mathbf{e}^{G_{2}(t_{k+1}^{i}-\tau)}d\tau\\
\end{split}
\end{equation}
which is not  equivalent to \eqref{Au2} and may  have the poorer performance than \eqref{Au2} since $u_{i}(t)$ in \eqref{Au3} depends on the sampled state $\phi_{i}(t_{k}^{i})$  instead of the continuous state $\phi_{i}(t)$.
In summary,  the key feature of our control law is that what we have designed is the same as what we will implement in the digital platform.
\end{Remark}
\begin{Remark}

It is worth mentioning  that  a special
case of this paper with $N = 1$ was studied in \cite{LiuHuang2017a}.  What makes the  current paper interesting is that
both our control law and event-triggered mechanism must be distributed in the sense that
the control law and the event-triggered mechanism of each subsystem can only make use of the
information of itself and its neighbors. This constraint poses significant difficulty in seeking
for a suitable control law and a suitable event-triggered mechanism. Moreover, in order to
obtain our result under the topological assumption that the digraph of the system is
connected and directed, our closed-loop system is a multi-input, multi-output coupled hybrid system. We
need to develop special skills to furnish a stability analysis of the closed-loop system.

\end{Remark}
\section{An Example}
Consider a class of Lorenz multi-agent systems taken from \cite{DongHuang2014a} as follows
\begin{equation}\label{Lorenzsystem2}
\begin{split}
 &\dot{z}_{1i}=c_{1i}z_{1i}-c_{1i}y_{i}\\
 &\dot{z}_{2i}=c_{2i}z_{2i}+z_{1i}y_{i}\\
 &\dot{y}_{i}= c_{3i}z_{1i}-y_{i}-z_{1i}z_{2i} + b_{i}u_{i}\\
 &e_{i}=y_{i}-v_{1},~ i=1,2,3,4
 \end{split}
\end{equation}
where  $c_{i}\triangleq\mbox{col} (c_{1i}, c_{2i},c_{3i},b_{i})$ is a constant parameter vector satisfying $c_{1i}<0$, $c_{2i}<0$ and $b_{i}>0$ for $i=1,2,3,4$.
The leader system takes the form of (\ref{exosystem1}) with $S=\left[
  \begin{array}{cc}
    0 & 1\\
    -1 & 0 \\
  \end{array}\right]$.  Clearly, Assumption \ref{Ass2.1} is satisfied.
The uncertain parameter $c_{i}$ is expressed as $c_{i}=\bar{c}_{i}+w_{i}$, where
$\bar{c}_{i}=\mbox{col}(\bar{c}_{1i},\bar{c}_{2i},\bar{c}_{3i},\bar{b}_{i})=\mbox{col}(-6,-8,1,2)$
is the nominal value of $c_{i}$, and $w_{i}=\mbox{col}(w_{1i},w_{2i},w_{3i},w_{4i})$
is the uncertainty of $c_{i}$ for $i=1,2,3,4$. Let $w=\mbox{col}(w_{1},w_{2},w_{3},w_{4})$. We assume that
$w\in\mathbb{W}=\{w\;|\;w\in\mathbb{R}^{16}, |w_{ji}|\leq1,~i,j=1,2,3,4\}$, and $v\in\mathbb{V}=\{v\;|\;v\in\mathbb{R}^{2}, |v_{i}|\leq1,~i=1,2\}$.

\begin{figure}[H]
  \centering
    \includegraphics[width=1.5in]{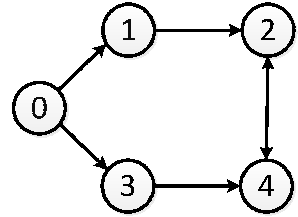}
  \caption{Communication network topology}
  \label{g} 
\end{figure}

The communication network topology is described  in Figure \ref{g}, where node $0$ is associated with the leader and the other nodes are associated with the followers.
It is easy to see that Assumption \ref{Ass3.2} is satisfied.

 Like in \cite{DongHuang2014a},  we have
\begin{equation*}\label{regulator2}
\begin{split}
&\textbf{y}_{i}(v,w)=v_{1},~\textbf{z}_{1i}(v,w)=r_{11i}(w)v_{1}+r_{12i}(w)v_{2}\\
&\textbf{z}_{2i}(v,w)=r_{21i}(w)v_{1}^{2}+r_{22i}(w)v_{2}^{2}+r_{23i}(w)v_{1}v_{2}\\
&\textbf{u}_{i}(v,w)=r_{31i}(w)v_{1}\!+\!r_{32i}(w)v_{2}\!+\!r_{33i}(w)v_{1}^{3}\!+\!r_{34i}(w)v_{2}^{3}\\
&~~~~~~~~~~~~~+r_{35i}(w)v_{1}^{2}v_{2}+r_{36i}(w)v_{1}v_{2}^{2},~i=1,2,3,4
 \end{split}
\end{equation*}
where the coefficients  can be found in \cite{DongHuang2014a}. Clearly, Assumptions \ref{Ass2.3} and \ref{Ass2.4} are satisfied.
Also, we can further verify that, for $i=1,2,3,4$,
\begin{equation}\label{Ui3}
\begin{split}
\dfrac{d^{4}\textbf{u}_{i}(v,w)}{dt^{4}}+10\dfrac{d^{2}\textbf{u}_{i}(v,w)}{dt^{2}}+9\textbf{u}_{i}(v,w)=0.
 \end{split}
\end{equation}
Thus, for $i=1,2,3,4$, we have
\begin{equation}\label{taui3}
\begin{split}
&\tau_{i}(v,w)=\mbox{col} (\textbf{u}_{i},\dot{\textbf{u}}_{i},\textbf{u}_{i}^{(2)},\textbf{u}_{i}^{(3)}),\\
&\Phi_{i}=\left[
                                                                                                              \begin{array}{cccc}
0 & 1 & 0 & 0 \\
0 & 0 & 1 & 0 \\
0& 0 & 0 & 1 \\
-9& 0 & -10 & 0 \\
                                                                                                              \end{array}
                                                                                                            \right]
,\ \ \Gamma_{i}=\left[
                  \begin{array}{c}
                    1 \\
                    0 \\
                    0 \\
                    0 \\
                  \end{array}
                \right]^{T}.
 \end{split}
\end{equation}
Choose the controllable pair $(M_{i},Q_{i})$  as follows,
\begin{equation*}
\begin{split}
M_{i}=\left[
                                                                                                              \begin{array}{cccc}
0 & 1 & 0 & 0 \\
0 & 0 & 1 & 0 \\
0 & 0 & 0 & 1 \\
-4 & -12 & -13 & -6 \\
                                                                                                              \end{array}
                                                                                                            \right],\
                                                                                                            Q_{i}=\left[
                  \begin{array}{c}
                    0 \\
                    0 \\
                    0 \\
                    1 \\
                  \end{array}
                \right].
\end{split}
\end{equation*}
By solving the Sylvester equation \eqref{Sylvester}, 
 we have $\Psi_{i}=\Gamma_{i}T^{-1}_{i}=[-5, 12, 3, 6]$. Then we perform the coordinate transformation (\ref{transformation}) and get the following augmented system,
\begin{equation}\label{Lorenzsystem3}
\begin{split}
 &\dot{\bar{z}}_{i}=\bar{f}_{i}(\bar{z}_{i},e_{i},\mu)\\
 &\dot{\bar{\eta}}_{i}=M_{i}\bar{\eta}_{i}+M_{i}Q_{i}b_{i}^{-1}e_{i}-Q_{i}b_{i}^{-1}\bar{g}_{i}(\bar{z}_{i},e_{i},\mu)\\
 &\dot{e}_{i}=\bar{g}_{i}(\bar{z}_{i},e_{i},\mu)\!+\!b_{i}\Psi_{i}\bar{\eta}_{i}\!+\!\Psi_{i}Q_{i}e_{i}\!+\!b_{i}\bar{u}_{i}\!+\!b_{i}\Psi_{i}\tilde{\eta}_{i}\\
 \end{split}
\end{equation}
where $\bar{z}_{i}=\mbox{col}(\bar{z}_{1i},\bar{z}_{2i})$
\begin{equation*}
\begin{split}
&\bar{f}_{i}(\bar{z}_{i},e_{i},\mu)=\left[
                                     \begin{array}{c}
                                       c_{1i}\bar{z}_{1i}-c_{1i}e_{i} \\
                                       c_{2i}\bar{z}_{2i}+(\bar{z}_{1i}+\textbf{z}_{1i})(e_{i}+v_{1})-\textbf{z}_{1i}v_{1} \\
                                     \end{array}
                                   \right]\\
&\bar{g}_{i}(\bar{z}_{i},e_{i},\mu)=c_{3i}\bar{z}_{1i}-e_{i}-\bar{z}_{1i}\bar{z}_{2i}-\textbf{z}_{1i}\bar{z}_{2i}-\bar{z}_{1i}\textbf{z}_{2i}.\\
 \end{split}
\end{equation*}
For the $\bar{z}_{i}$-subsystem, choose the Lyapunov function candidate as follows:
\begin{equation}
  \begin{aligned}
V_{1i}(\bar{z}_{i})=\frac{\hbar_{i}}{2}\bar{z}_{1i}^{2}+\frac{\hbar_{i}}{4}\bar{z}_{1i}^{4}+\frac{1}{2}\bar{z}_{2i}^{2}
  \end{aligned}
\end{equation}
for some sufficiently large $\hbar_{i}>0$. It is possible to show that, for all $\mu\in\Omega$,  all $\bar{z}_{i}$, and  all $e_{i}$ 
\begin{equation}
  \begin{aligned}
\frac{\partial V_{1i}(\bar{z}_{i})}{\partial \bar{z}_{i}}\bar{f}_{i}(\bar{z}_{i},e_{i},\mu)\leq&-\ell_{1i}\bar{z}_{1i}^{2}-\ell_{2i}\bar{z}_{1i}^{4}-\ell_{3i}\bar{z}_{2i}^{2}\\
&+\ell_{4i}e_{i}^{2}+\ell_{5i}e_{i}^{4}\\
  \end{aligned}
\end{equation}
for some constants $\ell_{si}>0$, $s=1,\cdots,5$, $i=1,\cdots,4$. That is to say, Assumption \ref{Ass3.1} is also satisfied.

Thus, by Theorem \ref{Theorem1}, we can design a distributed output feedback  control law of the form \eqref{ui3} with $\rho_{i}(e_{vi}(t_{k}^{i}))=10(e_{vi}^{6}(t_{k}^{i})+1)$, and a distributed output-based event-triggered mechanism of the form \eqref{trigger2} with $\sigma=0.01$, and $\delta=0.02$ or $0.002$.



\begin{figure}[H]
\centering
\includegraphics[scale=0.53]{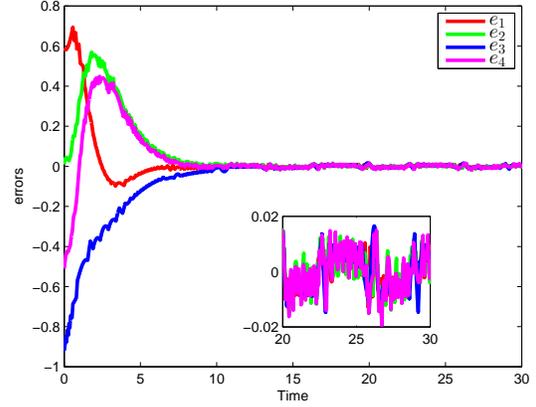}
\caption{Tracking errors of all followers for $\delta=0.02$ } \label{error1}
\end{figure}
\begin{figure}[H]
\centering
\includegraphics[scale=0.53]{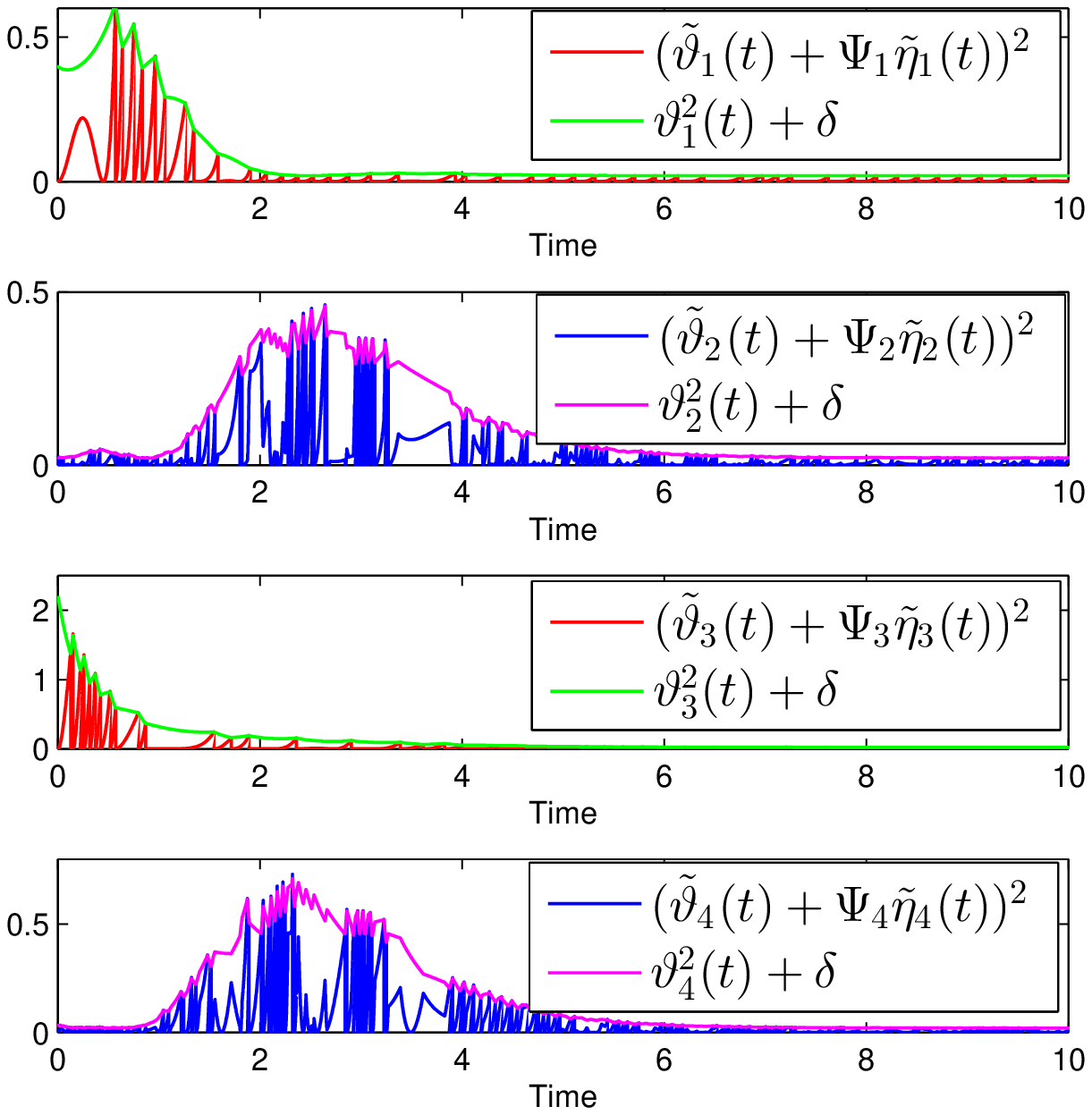}
\caption{Event-triggered  conditions  for $\delta=0.02$ } \label{condition1}
\end{figure}
\begin{figure}[H]
\centering
\includegraphics[scale=0.53]{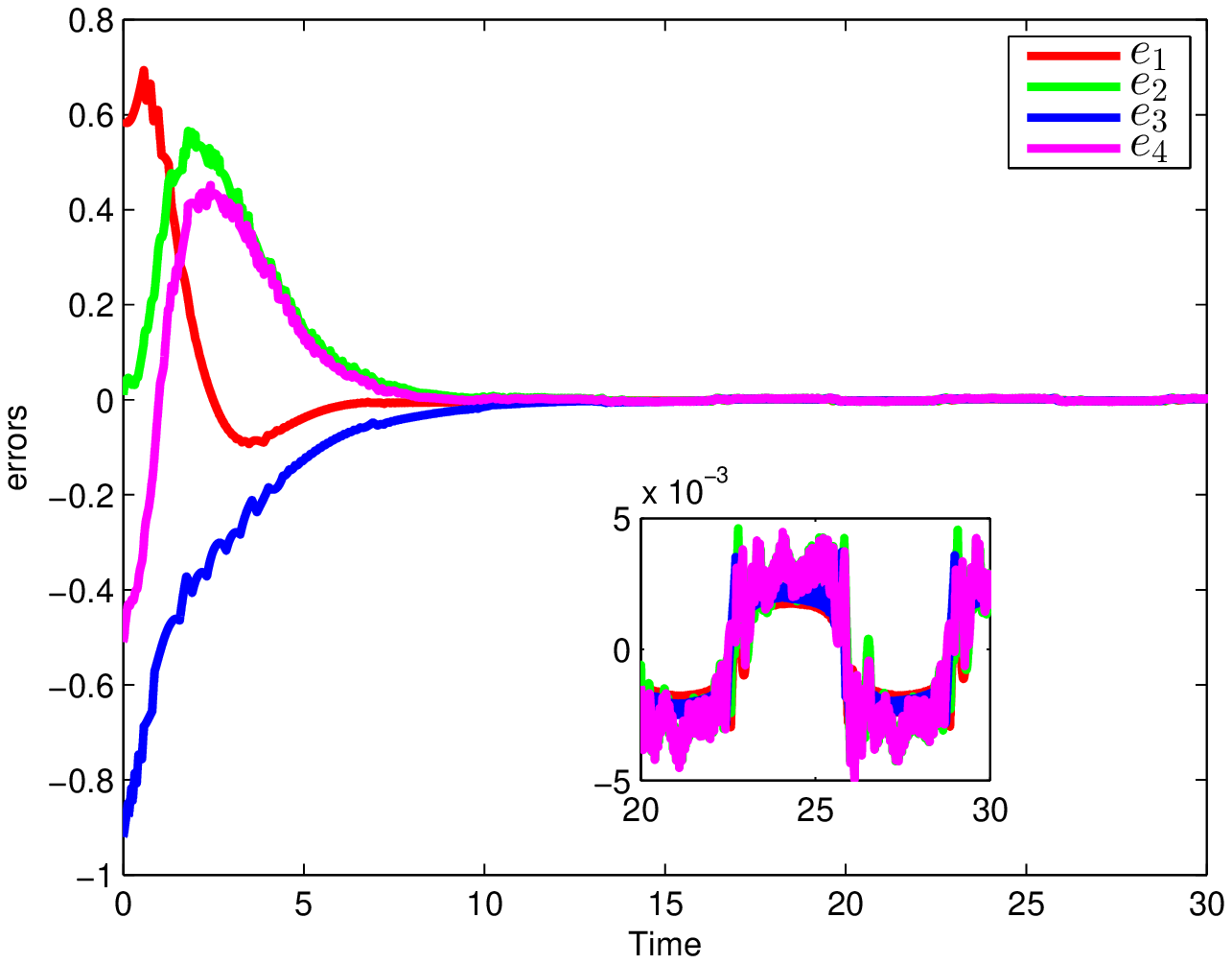}
\caption{Tracking errors of all followers  for $\delta=0.002$ } \label{error2}
\end{figure}
\begin{figure}[H]
\centering
\includegraphics[scale=0.53]{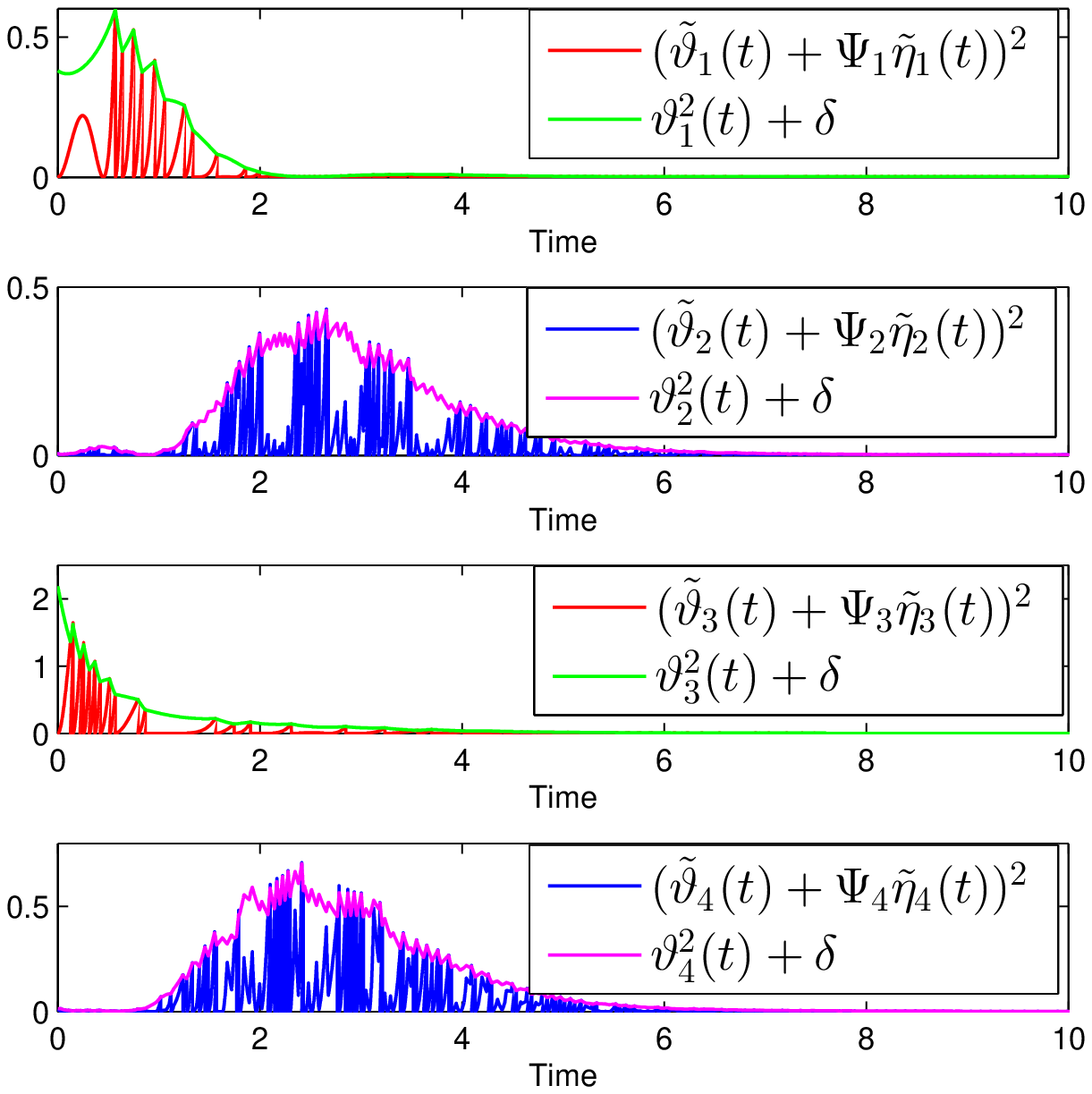}
\vspace{-0.5cm}
\caption{Event-triggered  conditions  for $\delta=0.002$ } \label{condition2}
\end{figure}
\begin{table}[H]
  \begin{center}
    \caption{Event-triggered numbers of all agents.}\label{table}
  \scalebox{0.78}{
    \begin{tabular}{|c|c|c|c|c|c|c|} \hline
      \multirow{2}{*}{Design parameters} & \multirow{2}{*}{Time} & \multicolumn{4}{c|}{Triggering numbers for each agent} \\
      \cline{3-6} &  & Agent 1 & Agent 2 & Agent 3 & Agent 4\\ \hline
      \multirow{2}{*}{$\sigma=0.01$, $\delta=0.02$} & \multirow{2}{*}{0-10s}   & \multirow{2}{*}{44} & \multirow{2}{*}{100} & \multirow{2}{*}{44} & \multirow{2}{*}{117} \\
       &    &  &  &  &  \\ \hline
      \multirow{2}{*}{$\sigma=0.01$, $\delta=0.002$} & \multirow{2}{*}{0-10s}   & \multirow{2}{*}{93} & \multirow{2}{*}{128} & \multirow{2}{*}{73} & \multirow{2}{*}{154} \\
       &    &  &  &  &  \\ \hline
    \end{tabular}}
  \end{center}
\end{table}

Simulation is performed with
\begin{equation*}
\begin{split}
 &w_{1}=[0.4,-0.7,0.6,-0.2]^{T}\\
 &w_{2}=[-0.3,0.2,0.5,0.4]^{T}\\
 &w_{3}=[-0.6,-0.4,0.8,-0.5]^{T}\\
 &w_{4}=[0.3,-0.4,0.6,0.3]^{T}\\
 \end{split}
\end{equation*}
 and the following initial conditions
\begin{equation*}
\begin{split}
 &v(0)=[0.88, -0.48]^{T}\\
 &[z_{11}(0),z_{21}(0),y_{1}(0)]=[1.38, 0.22, 1.47]\\
 &[z_{12}(0),z_{22}(0),y_{2}(0)]=[0.75, 0.04, 0.89]\\
 &[z_{13}(0),z_{23}(0),y_{3}(0)]=[0.01, -0.01, -0.04]\\
 &[z_{14}(0),z_{24}(0),y_{4}(0)]=[0.24, 0.01, 0.37]\\
 &\eta_{1}(0)=[ -0.75, 0.94, -0.14, -1.51]^{T}\\
 &\eta_{2}(0)=[-0.21, -1.16, 0.57, 1.78]^{T}\\
 &\eta_{3}(0)=[0.83, -0.35, 0.52, -1.33]^{T}\\
 &\eta_{4}(0)=[-0.50, -1.27,  0.83, 1.53]^{T}.\\
 \end{split}
\end{equation*}

Table \ref{table} shows the event-triggered numbers of all agents for $\delta=0.02$ and $\delta=0.002$. Figures \ref{error1} and \ref{error2} show the tracking errors of all followers for $\delta=0.02$ and $\delta=0.002$. Figures \ref{condition1} and \ref{condition2} show the event-triggered conditions of all followers for $\delta=0.02$ and $\delta=0.002$. These simulation results confirm that $\lim_{t\rightarrow\infty}\sup |e_{i}(t)|\leq 0.02 $ for  $\delta=0.02$ and $\lim_{t\rightarrow\infty} |e_{i}(t)|\leq 5\times10^{-3}$ for $\delta=0.002$. Moreover, the event-triggered numbers of all agents for $\delta=0.02$ are less than those for $\delta=0.002$.

\section{Conclusion}
In this paper, we have studied the event-triggered cooperative global robust practical output regulation problem for a class of nonlinear multi-agent systems by a distributed output-based event-triggered control law together with a distributed output-based event-triggered mechanism, and proved that the Zeno behavior can be avoided. As our result applies to  nonlinear multi-agent systems with both external disturbances and unknown parameters and
does not require the nonlinear functions satisfy the global Lipchitz condition, it significantly enlarges the systems that can be handled by the existing approaches. As a natural extension of the result in this paper, we are considering the same problem under switching network topologies.

\end{document}